\documentclass[12pt,a4paper]{article}

\usepackage[font=footnotesize,width=0.8\textwidth]{caption}

\usepackage{algorithm}
\usepackage{algorithmic}
\usepackage{amsmath,amssymb}
\usepackage{amsfonts}
\usepackage{amsbsy}
\usepackage{amsthm} 
\usepackage{array}
\usepackage[english]{babel}
\usepackage{bibentry}
\usepackage{bm}
\usepackage{booktabs} 
\usepackage{calc}
\usepackage{cite}
\usepackage{color}
\usepackage{fancyhdr}
\usepackage{graphicx}
\usepackage{graphics}
\usepackage[bookmarks=true,pdftex=false,bookmarksopen=true,pdfborder={0,0,0}]{hyperref}
\usepackage{hyphenat} 
\usepackage[latin1]{inputenc}
\usepackage{latexsym}
\usepackage{listings}
\usepackage{longtable}
\usepackage{lscape}
\usepackage{mathtools}
\usepackage{mdwmath}
\usepackage{mdwtab}
\usepackage{multirow}
\usepackage[numbers]{natbib}
\usepackage{paralist}
\usepackage{pgf}
\usepackage{placeins}
\usepackage{psfrag}
\usepackage{rotating}
\usepackage{stackengine}
\usepackage{stmaryrd}
\usepackage{subcaption} 
\usepackage{tabularx}
\usepackage{tikz}
\usepackage{url}
\usepackage{verbatim}

\usetikzlibrary{matrix,arrows}  
\usetikzlibrary{shapes}         

\newtheorem{definition}{Definition}[section]

\newtheorem{lemma}{Lemma}[section]
\newtheorem{theorem}[lemma]{Theorem}
\newtheorem{corollary}[lemma]{Corollary}

\newenvironment{proof*}[1][Proof.]{\begin{trivlist}
\item[\hskip \labelsep {\bfseries #1}]}{\end{trivlist}}

\newenvironment{remark*}[1][Remark.]{\begin{trivlist}
\item[\hskip \labelsep {\itshape #1}]}{\end{trivlist}}

\newcommand{\tb}{\textbf}
\let\dd\partial            

\newcommand{\de}{\mathop{}\!\mathrm{d}} 

\AtBeginDocument{
  \let\div\relax
  \DeclareMathOperator{\div}{div}  
}

\AtBeginDocument{
  
}

\AtBeginDocument{
  
}

\AtBeginDocument{
  \let\sup\relax
  \DeclareMathOperator*{\sup}{sup}  
}

\AtBeginDocument{
  
}

\AtBeginDocument{
  
}

\newcommand{\closure}[2][3]{
{}\mkern#1mu\overline{\mkern-#1mu#2}} 

\newcommand{\ol}[1]{%
\mkern 1.5mu\overline{\mkern-1.5mu#1\mkern-1.5mu}\mkern 1.5mu} 

\AtBeginDocument{
  \let\diam\relax
  \DeclareMathOperator{\diam}{diam}  
}

\newcommand\rst[2]{{%
  \left.\kern-\nulldelimiterspace 
  #1 
  \vphantom{\big|} 
  \right|_{#2} 
  }}

\stackMath
\newcommand\stackti[2][1]{%
 \def\useanchorwidth{T}%
  \ifnum#1>1%
    \stackon[0pt]{\stackti[\numexpr#1-1\relax]{#2}}{\scriptscriptstyle\sim}%
  \else%
    \stackon[1pt]{#2}{\scriptscriptstyle\sim}%
  \fi%
}

\renewcommand{\texttt}[1]{%
  \begingroup
  \ttfamily
  \begingroup\lccode`~=`/\lowercase{\endgroup\def~}{/\discretionary{}{}{}}%
  \begingroup\lccode`~=`[\lowercase{\endgroup\def~}{[\discretionary{}{}{}}%
  \begingroup\lccode`~=`.\lowercase{\endgroup\def~}{.\discretionary{}{}{}}%
  \catcode`/=\active\catcode`[=\active\catcode`.=\active
  \scantokens{#1\noexpand}%
  \endgroup
}

\begin{document}

\title{Convergence analysis of a cell centered\\ finite volume diffusion operator\\ on non-orthogonal polyhedral
  meshes}

\author{Luca Bonaventura$^{(1)}$,\ \ Alessandro Della Rocca$^{(1),(2)}$} \maketitle

\begin{center}
{\small
$^{(1)}$ MOX -- Modelling and Scientific Computing, \\
Dipartimento di Matematica, Politecnico di Milano \\
Via Bonardi 9, 20133 Milano, Italy\\
{\tt luca.bonaventura@polimi.it}\\
{\tt alessandro.dellarocca@polimi.it}
}
\end{center}

\begin{center}
{\small $^{(2)}$ Tenova S.p.A., \\ Global R\&D, \\ Via Albareto 31, 16153 Genova, Italy\\ {\tt
    alessandro.dellarocca@tenova.com} }
\end{center}

\date{}

\noindent
{\bf Keywords}:  Finite volume, Cell centered methods, Convergence, Diffusion, Unstructured  meshes 

\vspace*{1.0cm}

\noindent
{\bf AMS Subject Classification}: 65M08, 65N08, 65N12, 65Z05, 76R50

\vspace*{1.0cm}

\pagebreak

\abstract{A simple but successful strategy for building a discrete diffusion operator in finite volume schemes
  of industrial use is to correct the standard two-point flux approximation with a term accounting for the
  local mesh non-orthogonality. Practical experience with a variety of different mesh typologies, including
  non-orthogonal tetrahedral, hexahedral and polyhedral meshes, has shown that this discrete diffusion
  operator is accurate and robust whenever the mesh is not too distorted and sufficiently regular. In this
  work, we show that this approach can be interpreted as equivalent to introducing an anisotropic operator
  that accounts for the preferential directions induced by the local mesh non-orthogonality. This allows to
  derive a convergence analysis of the corrected method under a quite weak global assumption on mesh
  distortion.  This convergence proof, which is obtained for the first time for this finite volume method
  widely employed in industrial applications, provides a reference framework on how to interpret some of its
  variants commonly implemented in commercial finite volume codes.  Numerical experiments are presented that
  confirm the accuracy and robustness of the results. Furthermore,  we also show empirically that a least square approach to the gradient computation  can
  provide second order convergence even when the mild   non-orthogonality condition on the mesh is violated.
}

\pagebreak


\section{Introduction}
 \label{sec:introduction} \indent

Finite volume methods have been extremely popular in computational fluid dynamics (CFD) in the past and they
still are an area of active research in numerical mathematics. Among the many different developments in this
field, we recall the finite volume element scheme \cite{cai:1990}, the multi-point flux approximation schemes
(MPFA) \cite{aavatsmark:1998a}, \cite{aavatsmark:1998b}, \cite{aavatsmark:2008}, or more recent variants, such
as the mixed finite volume scheme (MFV) \cite{droniou:2006}, the hybrid finite volume scheme (HFV)
\cite{eymard:2007c}, \cite{eymard:2010} and the discrete duality finite volume schemes (DDFV)
\cite{andreianov:2013},\cite{andreianov:2012}, \cite{coudiere:2011}, \cite{coudiere:2010},
\cite{hermeline:2007}, \cite{hermeline:2009}, \cite{krell:2012}.  All these methods share indeed many common
features, as discussed in \cite{droniou:2010}.

In this work, however, we will focus on cell centered schemes, in which a single unknown is associated to each
mesh cell.  Cell centered finite volume methods are widely employed in industrial codes \cite{Fluent},
\cite{OpenFOAM} for a number of practical reasons. Indeed, they rely on relatively simple data structures,
even for general unstructured meshes, and they allow for easy treatment of boundary conditions at
singular boundary points, such as inner corners, while effectively handling general shapes of the
computational domain.  Cell centered methods can be naturally parallelized by domain decomposition techniques,
guaranteeing minimal interprocessor communications, especially in their low order variants, due to the use of
discrete operators built from local stencils. They allow an easy implementation of locally adaptive multilevel
refinement strategies and they can be easily equipped with very efficient geometric multigrid procedures
\cite{wesseling:2001}. Finally, cell centered finite volume methods also allow an immediate extension to
nonlinear coupled problems \cite{eymard:2009b}.

Known drawbacks of cell centered schemes are the reduced accuracy in strongly heterogeneous diffusion problems
\cite{eymard:2010} with respect to MPFA, MFV or HFV schemes, as well as the only asymptotical recovery of the
discrete Stokes formula, in contrast with the exact discrete property provided for example by DDFV schemes
\cite{coudiere:2011}. On the other hand, MPFA, MFV, HFV and DDFV schemes achieve such properties by
introducing additional unknowns at selected mesh locations, thus implying an additional cost with respect to
cell centered discretizations. It is still an open question if similar accuracy improvements can be obtained
from cell centered schemes by introducing additional unknowns through local mesh refinement.

For these reasons, it is important to understand the analytical behaviour of cell centered finite volume
discretizations on the typical non-orthogonal meshes practicaly required for industrial applications
\cite{Fluent}, \cite{OpenFOAM}. For these applications, the so-called \emph{Gauss corrected} scheme, widely
adopted by finite volume practitioners \cite{ferziger:2002}, \cite{mathur:1997}, \cite{moukalled:2016},
\cite{muzaferija:1997}, \cite{tsui:2006}, appears to be a simple, robust and sufficiently accurate
option. Notice that this scheme can also be interpreted as a specific realization of the recently introduced
asymmetric gradient discretization method \cite{droniou:2017}.

To the best of the authors' knowledge, the convergence properties of this finite volume method have never been
analyzed in the case of non-orthogonal meshes.  Indeed, convergence analyses of finite volume schemes for
diffusion operators on unstructured mesh types are usually limited to polyhedral meshes satisfying an
orthogonality condition \cite{eymard:2000}, \cite{eymard:2006b}. This is quite restrictive in practice, since
none of the robust mesh generators usually adopted for pre-processing of industrial configurations are able
to guarantee this condition.

In this work, we show that it is possible to prove the convergence of the \emph{Gauss corrected} scheme on
unstructured meshes satisfying a global and rather weak mesh regularity condition. 
 This goal is achieved adapting the approach used in \cite{eymard:2006b} for the
convergence analysis of a cell-centered finite volume scheme for anisotropic diffusion problems on orthogonal
meshes. A preliminary version of these results has been presented in
\cite{dellarocca:2018}. Furthermore,  we also show empirically that a least square approach to the gradient computation  can  provide second order convergence even when the mild   non-orthogonality condition on the mesh is violated.
It is to be remarked that existing convergence proofs for finite volume methods on
non-orthogonal meshes either involve discretization schemes not guaranteeing local flux conservativity
\cite{eymard:2009b}, \cite{eymard:2010}, or DDFV schemes employing additional degrees of freedom
\cite{andreianov:2013}, \cite{andreianov:2012}, or two-dimensional diamond schemes on meshes satisfying more
restrictive regularity conditions \cite{coudiere:1999}.
\noindent
We will focus here on the isotropic steady state diffusion equation
\begin{subequations}
 \begin{align}
  - \div(\alpha \nabla \underline{u}) &= f, \quad \text{in}\; \Omega,\\
  \underline{u} &= 0, \quad \text{on}\; \dd \Omega.
 \end{align}
 \label{eq:isotropic_diffusion_problem}%
\end{subequations}
We will assume that $\alpha : \Omega \rightarrow \mathbb{R}$ is a measurable function, $\alpha \in
L^\infty(\Omega)$,   such that $0 < \alpha_0 \leq \alpha(\bm{x})$ for
a.e. $\bm{x} \in \mathbb{R}^d$, with $\alpha_0 \in \mathbb{R}$, and for $f \in L^2(\Omega)$.
The classical weak formulation of problem \eqref{eq:isotropic_diffusion_problem} consists in finding
$\underline{u} \in H_0^1(\Omega) $ such that
\begin{equation}
  \int_{\Omega} \alpha(\bm{x}) \, \nabla \underline{u}(\bm{x}) \cdot \nabla v(\bm{x}) \,\de\bm{x}
    = \int_{\Omega} f(\bm{x}) \, v(\bm{x}) \,\de\bm{x},
   \quad \forall v \in H_0^1(\Omega).
 \label{eq:isotropic_diffusion_weak_problem}
\end{equation}
Rather than proving convergence directly for the finite volume scheme associated to the strong problem
formulation \eqref{eq:isotropic_diffusion_problem}, we will identify a discrete weak formulation underlying
the finite volume scheme and then to prove convergence of its solution to that of the continuous weak problem
\eqref{eq:isotropic_diffusion_weak_problem}.

The paper is organized as follows. In section \ref{sec:meshes_spaces}, several fundamental definitions of mesh
related quantities and discrete functional spaces are introduced.  In section \ref{sec:fvm_diffusion}, the
cell centered finite volume method that is the focus of our analysis is presented. In section
\ref{sec:underlying_weak_form}, the discrete weak formulation is recovered and in section
\ref{sec:convergence_analysis}, the convergence analysis of the \emph{Gauss corrected} scheme is presented. In
section \ref{sec:numerical_results}, the results of some numerical experiments are reported.  A proposal
to overcome the constraints on the mesh for some specific three-dimensional mesh types is introduced in
section \ref{sec:least_squares_tests}. Finally, in section \ref{sec:conclusions} some conclusions are drawn
and some future developments are outlined.

\section{Meshes and discrete spaces}
 \label{sec:meshes_spaces} \indent

The finite volume method is a mesh-based discretization technique suitable for any number of space dimensions,
but in this work we only consider the $d=3$ case.  Since the computational domains of practical interest are
usually of complex geometry, the focus here is on meshes composed of arbitrarily shaped polyhedral cells, in
the sense of the formal definition below, see also \cite{droniou:2016a}, \cite{eymard:2010}.
\begin{definition} \textbf{(Polyhedral mesh)}:
Let $\Omega$ be a bounded, open polyhedral subset of $\mathbb{R}^d. $   A polyhedral mesh for
$\Omega$ is denoted by $\mathcal{D}=(\mathcal{M},\mathcal{F},\mathcal{P},\mathcal{V})$, where the quadruple
includes:
\begin{enumerate}
\item $\mathcal{M}$ is a finite family of non-empty, connected, polyhedral, open,
 disjoint subsets of $\Omega$ called cells (or control volumes), such that $\,\overline{\Omega}
 = \cup_{K \in \mathcal{M}} \closure{K}$. For any $K \in \mathcal{M}$, $\partial K = \closure{K} \setminus K$
 is the boundary of $K$, $|K| > 0$ denotes the measure of $K$, and $h_K = \diam(K)$ is the diameter of $K$,
 that is the maximum distance between two points in $K$.
\item $\mathcal{F} = \mathcal{F}_{int} \cup \mathcal{F}_{ext}$ is a finite family of disjoint
 subsets of $\,\overline{\Omega}$ representing the faces. Let $\mathcal{F}_{int}$ be the set of interior faces
 such that, for all $\sigma \in \mathcal{F}_{int}$, $\sigma$ is a non-empty open subset of a hyperplane in
 $\mathbb{R}^d$ with $\sigma \subset \Omega$, and let $\mathcal{F}_{ext}$ be the set of boundary faces such
 that, for all $\sigma \in \mathcal{F}_{ext}$, $\sigma$ is a non-empty open subset of $\partial \Omega$.  It
 is assumed that, for any $K \in \mathcal{M}$, there exists a subset $\mathcal{F}_K \subset \mathcal{F}$ such
 that $\partial K = \cup_{\sigma \in \mathcal{F}_K} \overline{\sigma}$. The set of cells sharing one face
 $\sigma$ is $\mathcal{M}_{\sigma} = \{ K \in \mathcal{M} : \sigma \in \mathcal{F}_K\}$.  It is assumed that,
 for all $\sigma \subset \mathcal{F}$, either $\mathcal{M}_{\sigma}$ has exactly two elements and then
 $\sigma \subset \mathcal{F}_{int}$, or $\mathcal{M}_K$ has exactly one element and then
 $\sigma \subset \mathcal{F}_{ext}$.  For all $\sigma \in \mathcal{F}$, $|\sigma|>0$ denotes the
 $(d-1)$-dimensional measure of $\sigma$, and $\overline{\tb{x}}_{\sigma}$ is the barycenter of $\sigma$.
\item $\mathcal{P} = (\tb{x}_K)_{K \in \mathcal{M}}$ is a family of points of $\Omega$ indexed by $\mathcal{M}$,
 such that for all $K \in \mathcal{M}$, $\bm{x}_K \in K$ and it is called the center of $K$, possibly
 corresponding to its barycenter. It is assumed that all cells $K \in \mathcal{M}$ are $\bm{x}_K$-star-shaped,
 in the sense that if $\bm{x} \in K$, then the line segment $[\bm{x}_K, \bm{x}] \subset K$.
\item $\mathcal{V}$ is the finite set of vertices of the mesh.
 For $K \in \mathcal{M}$, $\mathcal{V}_K$ collects all the vertices belonging to
 $\closure{K}$, while for $\sigma \in \mathcal{F}$, $\mathcal{V}_{\sigma}$ collects all the vertices belonging
 to $\sigma$.
\end{enumerate}
The size of the polyhedral mesh is defined as $h_{\mathcal{D}} = \sup\{ h_K, K \in \mathcal{M} \}$.
 \label{def:polyhedral_mesh}
\end{definition}

\begin{figure}[t!]
  \centering
  \includegraphics[width=0.5\linewidth]{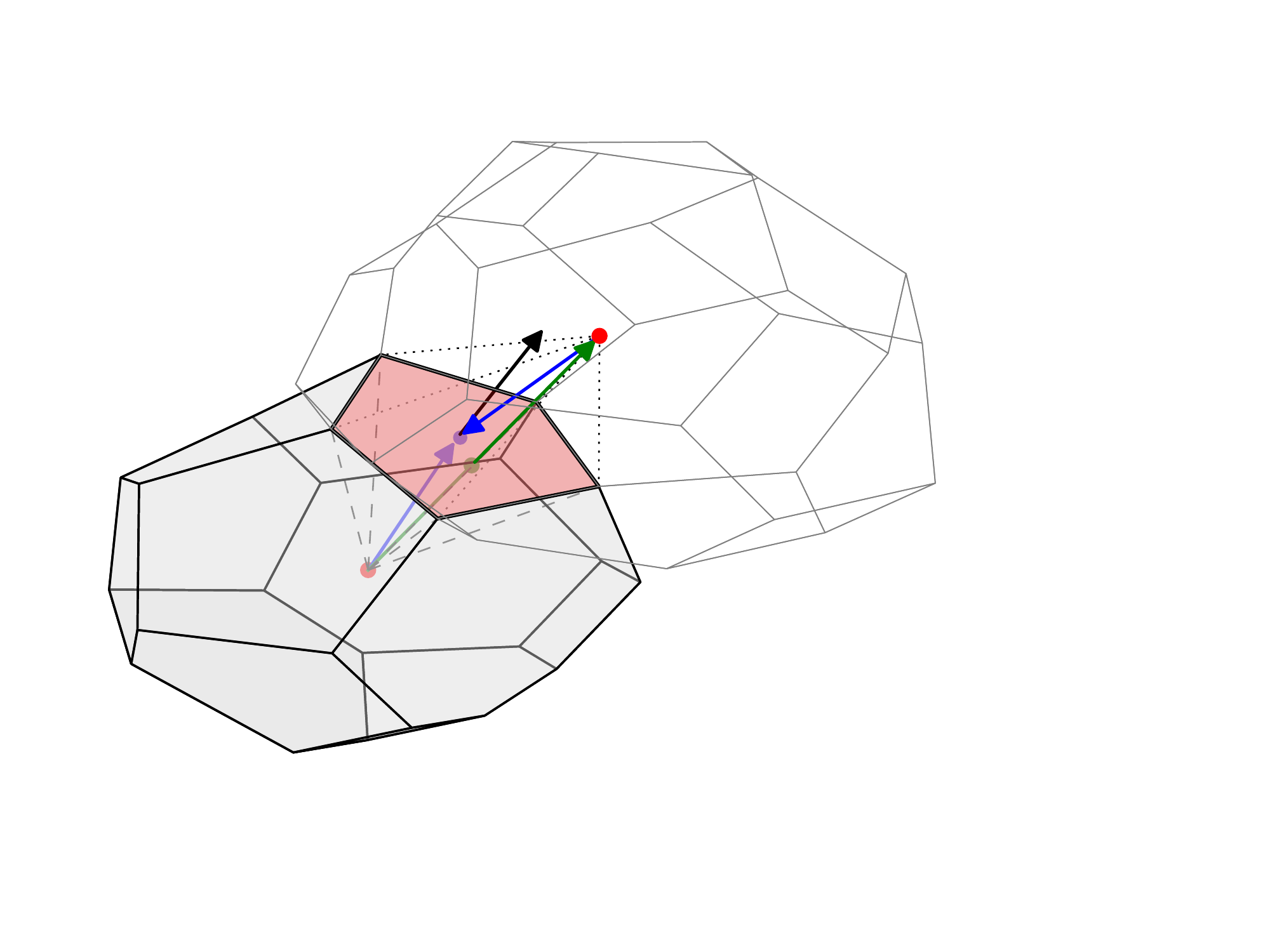}
  \caption{Non-orthogonal generally polyhedral mesh: the line segment $(\bm{x}_L - \bm{x}_K)$ is not aligned with
   the face normal unit vector $\bm{n}_{K,\sigma}$. Furthermore, the intersection point $\bm{y}_\sigma$ between
   the face $\sigma$ and the vector $(\bm{x}_L - \bm{x}_K)$ is not necessarily coincident with the face
   centroid $\bm{x}_\sigma$. The half diamonds $D_{K,\sigma}$ and $D_{L,\sigma}$ are  represented
   with dashed and dotted lines, respectively.}
 \label{fig:meshes}
\end{figure}

Furthermore, for any $K \in \mathcal{M}$ and for any $\sigma \in \mathcal{F}_K$, $\bm{n}_{K,\sigma}$ is the
constant unit vector normal to $\sigma$ and outward to $K$.  For any $K \in \mathcal{M}$, the set of
neighbors of $K$ is denoted by
\begin{equation}
  \mathcal{N}_K = \{ L \in \mathcal{M} \setminus \{K\},
                    \exists \sigma \in \mathcal{F}_{int}, \mathcal{M}_{\sigma}=\{K,L\} \}.
 \label{eq:cell_neighbours}
\end{equation}
Additionally $d_{K,\sigma}$ denotes the orthogonal distance between $\bm{x}_K$ and $\sigma \in \mathcal{F}_K$
\begin{equation}
  d_{K,\sigma} = (\bm{x} - \bm{x}_K) \cdot \bm{n}_{K,\sigma},  
 \label{eq:cell_center_face_distance}
\end{equation}
which is constant for all $\bm{x} \in \sigma. $ 
From the assumption that $K$ is $\bm{x}_K$-star-shaped,
it follows that   $d_{K,\sigma} > 0$  and  that it also holds:
\begin{equation}
  \sum_{\sigma \in \mathcal{F}_K} |\sigma| \,d_{K,\sigma} = d\, |K| \quad \forall K \in \mathcal{M}.
 \label{eq:cell_volume}
\end{equation}
For all $K \in \mathcal{M}$ and $\sigma \in \mathcal{F}_K$, $D_{K,\sigma}$ denotes the cone with vertex
$\bm{x}_K$ and basis $\sigma$, also called half-diamond, that is the volume defined by
\begin{equation}
  D_{K,\sigma} = \{ t\,\bm{x}_K + (1-t)\,\bm{y},\, t \in (0,1),\, \bm{y} \in \sigma  \}.
 \label{eq:face_cone}
\end{equation}
For all $\sigma \in \mathcal{F}$, $D_{\sigma} = \cup_{K \in \mathcal{M}_{\sigma}} D_{K,\sigma}$ denotes the
diamond associated to face $\sigma$, as in Figure \ref{fig:meshes}.

Definition \ref{def:polyhedral_mesh} covers a wide range of meshes, including meshes with non-convex cells,
with non-planar faces requiring triangulation, or with hanging nodes. Furthermore, Definition
\ref{def:polyhedral_mesh} also includes tetrahedral and hexahedral meshes as particular cases, as well as
meshes with wedge and pyramidal cells.

%

Finite volume methods are traditionally introduced in discrete functional spaces of piecewise constant
functions \cite{eymard:2000}. In recent analyses \cite{eymard:2007a}, \cite{eymard:2009a}, associated inner
products, norms and seminorms are exploited to recast the discrete flux balance equations into an equivalent
variational form, which naturally allows to derive stability estimates and to investigate the numerical
convergence of specific schemes \cite{eymard:2010}.
In the classical finite volume framework, the discrete flux balance equation corresponding to problem
\eqref{eq:isotropic_diffusion_problem} takes the form
\begin{equation}
  \sum_{\sigma \in \mathcal{F}_K} F_{K,\sigma}(u) = \int_{K} f(x) \,dx \quad \forall K \in \mathcal{M}
 \label{eq:anisotropic_diffusion_fvm}
\end{equation}
where the face flux is such that
\begin{equation*}
  F_{K,\sigma}\approx-\int_{\sigma} \alpha(x) \, \nabla u(x) \cdot \bm{n}_{K,\sigma} \, d\gamma(x)
\end{equation*}
and $d\gamma(x)$ denotes the infinitesimal face area element. A relevant feature of the
scheme is the flux conservativity property
\begin{equation}
   F_{K,\sigma}(u) + F_{L,\sigma}(u) = 0\,,
 \label{eq:flux_conservativity}
\end{equation}
which is assumed to hold for all interior faces $\sigma \in \mathcal{F}_{int}$, where $K$ and $L$ are the
cells sharing the face $\sigma$.

The convergence analysis of cell centered finite volume schemes on arbitrary polyhedral meshes 
\cite{eymard:2010} may also require to introduce the space $H_{\mathcal{D}}(\Omega) \subset
L^p(\Omega)$, which consists of the functions that are piecewise constant on each cell
$K \in \mathcal{M}$. For all $v \in H_{\mathcal{D}}(\Omega)$ and for all $K \in \mathcal{M}$, the constant
value of $v$ in $K$ is denoted by $v_K$. Consequently, discrete functional analysis results for the
convergence of finite volume schemes \cite{eymard:2000}, \cite{eymard:2010} can be exploited.

In addition, in order to introduce proper test functions to check the convergence of the discrete
solution to the continuous solution of the weak formulation, for all functions $\psi \in C(\Omega)$ a
projection operator $P_{\mathcal{D}} : C(\Omega) \rightarrow H_{\mathcal{D}}(\Omega)$ is defined, such that
$P_{\mathcal{D}} \psi = (\psi(\bm{x}_K))_{K \in \mathcal{M}}$.

\section{A cell centered diffusion scheme for non-orthogonal meshes}
 \label{sec:fvm_diffusion}

The vast majority of finite volume schemes for diffusive problems are based on the application of the discrete
Gauss theorem. The numerical approximation is derived as
\begin{equation}
  \int_{K} \nabla \!\cdot\! (\alpha \nabla u) \,\de\bm{x}
  = \sum_{\sigma \in \mathcal{F}_K} \int_{\sigma} \alpha\, \bm{n}_{K,\sigma} \!\cdot\! \nabla u \,\de\gamma
  \approx \sum_{\sigma \in \mathcal{F}_K} F_{K,\sigma}^{(d)}(u)
 \label{eq:discrete_gauss_diffusive_term}
\end{equation}
where the numerical flux through face $\sigma$ is computed as
\begin{equation}
  F_{K,\sigma}^{(d)}(u) = |\sigma| \,\alpha_{K,\sigma}\, \bm{n}_{K,\sigma} \cdot \nabla_{K,\sigma} u,
  \quad \forall \sigma \in \mathcal{F}_K
 \label{eq:discrete_gauss_diffusive_term}
\end{equation}
and depends on the definition of the the face normal gradient.  Usually, $\alpha_{K,\sigma}$ is approximated
by the surface value interpolation $I_\sigma \alpha$, obtained from standard interpolation schemes.  Linear
interpolation is often chosen to preserve second order accuracy, while harmonic interpolation is sometimes
selected, especially when the scalar diffusivity field $\alpha$ is strongly non-homogeneous
\cite{eymard:2000}.  A variety of alternative schemes can be constructed to approximate the face normal
gradient $\bm{n}_{K,\sigma} \cdot \nabla_{K,\sigma} u$, each with its own specific features. Most of them are
traditionally studied empirically, by directly testing them on specific meshes and representative flow
problems \cite{perezsegarra:2006b}, \cite{perezsegarra:2006a}.

The simplest scheme for the face normal gradient is represented by the two-point flux approximation
\cite{eymard:2000}
\begin{equation}
  \bm{n}_{K,\sigma} \!\cdot\! \nabla_{K,\sigma} u = \frac{u_L - u_K}{d_{K,\sigma} + d_{L,\sigma}}.
 \label{eq:face_gradient}
\end{equation}
Even though unconditionally monotone and coercive \cite{droniou:2014}, it is of limited accuracy on
unstructured meshes, where mesh non-orthogonality may lead to severe errors in the approximation of the
diffusion fluxes \cite{eigestad:2005}.  In order to compensate for the unavoidable non-orthogonality of
realistic unstructured meshes, a simple but effective solution is provided by the \emph{Gauss corrected}
scheme, which introduces a non-orthogonal correction term \cite{mathur:1997} in the two-point flux scheme,
thus obtaining the approximation
\begin{equation}
  \begin{aligned}
  \bm{n}_{K,\sigma} \!\cdot\! \nabla_{K,\sigma} u &= \frac{u_L - u_K}{d_{K,\sigma} + d_{L,\sigma}}\\
    &\hspace{1.0em} + \left( \bm{n}_{K,\sigma} - \frac{\bm{x}_L - \bm{x}_K}{\;\bm{n}_{K,\sigma} \cdot (\bm{x}_L - \bm{x}_K)} \right)
        \cdot \nabla_{\sigma} u.
  \end{aligned}
  \label{eq:face_gradient_corrected}
\end{equation}
Here, the first term corresponds to the two-point flux contribution in Eq.\eqref{eq:face_gradient}, expressing
the diffusion flux component in the direction of the line segment $(\bm{x}_L - \bm{x}_K)$, while the second
term accounts for the local mesh non-orthgonality across the face $\sigma$, expressed as the difference
between the correct face normal diffusive flux estimated from a proper face gradient and the diffusive flux
along the direction of $(\bm{x}_L - \bm{x}_K)$. In order to avoid oscillatory solutions
\cite{perezsegarra:2006a}, it is important that the gradient $\nabla_{\sigma}u$ at face $\sigma$ is evaluated
using a different scheme from the one employed in the first term of
Eq.\eqref{eq:face_gradient_corrected}. Thus, the gradient $\nabla_{\sigma} u$ is usually estimated at face
$\sigma$ by interpolation of the neighbouring cells gradients $\nabla_{K} u$ and $\nabla_{L} u$ for
$\mathcal{M}_{\sigma}=\{K,L\}$, which can be either the standard linear interpolation
%
%
or, for increased simplicity, the midpoint rule.
%
%
Indeed, if the cell derivatives are linear approximations, the diffusion flux will be more accurate than first
order on very regular meshes \cite{mathur:1997}.
The \emph{Gauss corrected} scheme allows more accurate approximations than the two-point flux
approximation \eqref{eq:face_gradient}, but it is not, in general, unconditionally coercive on
arbitrary unstructured meshes. As a consequence, on irregular meshes it may become a source of numerical
instability. On orthogonal grids, this scheme reduces to the classical two-point flux scheme, since the
correction term vanishes.

Following \cite{chenier:2009}, on a general unstructured polyhedral mesh like that of Definition
\ref{def:polyhedral_mesh}, the centered discrete gradient operator $\nabla_{\mathcal{D}} :
H_\mathcal{D}(\Omega) \rightarrow H_\mathcal{D}(\Omega)^d$ is defined as the piecewise constant function

\begin{equation}
  \nabla_K u = \frac{1}{|K|} \sum_{\sigma \in \mathcal{F}_K} |\sigma| \left( I_\sigma u - u_K \right) \bm{n}_{K,\sigma}
 \label{eq:discrete_gradient}
\end{equation}
for $ u \in H_\mathcal{D}(\Omega). $ Since  for any closed control volume the geometrical relations 
\begin{equation}
  \sum_{\sigma \in \mathcal{F}_K} |\sigma| \,\bm{n}_{K,\sigma} \cdot \bm{e}^{(i)}
  = \sum_{\sigma \in \mathcal{F}_K} |\sigma| \,n_{K,\sigma}^{(i)} = 0 \quad  i=1,\ldots,d
 \label{eq:closed_cell_relation:a}
\end{equation}
hold, with $\bm{n}_{K,\sigma} = n_{K,\sigma}^{(i)} \bm{e}^{(i)}$ , Eq.\eqref{eq:discrete_gradient} is also equal to
\begin{equation}
  \nabla_K u = \frac{1}{|K|} \sum_{\sigma \in \mathcal{F}_K} |\sigma|\, I_\sigma u \;\bm{n}_{K,\sigma},
 \label{eq:gauss_discrete_gradient}
\end{equation}
which is easily recognized as the finite volume discretization of the gradient based on the Gauss theorem
\cite{ferziger:2002}. For this reason, the gradient approximation in Eq.\eqref{eq:discrete_gradient} is often
identified as the Gauss gradient scheme.

The consistency of the discrete gradient in Eq.\eqref{eq:discrete_gradient} has been analyzed in
\cite{eymard:2010}. It stems directly  from the geometrical identity
\begin{equation}
  \sum_{\sigma \in \mathcal{F}_K} |\sigma| \bm{n}_{K,\sigma} (\bm{x}_\sigma - \bm{x}_K)^\intercal = |K| \bm{I},
    \quad \forall K \in \mathcal{M},
 \label{eq:cell_volume_tensor_identity}
\end{equation}
where $(\bm{x}_\sigma - \bm{x}_K)^\intercal$ is the transpose of the vector $(\bm{x}_\sigma - \bm{x}_K) \in
\mathbb{R}^d$, see Figure \ref{fig:meshes}, and $\bm{I} \in \mathbb{R}^d \times \mathbb{R}^d$ is the identity
matrix. For any affine function $\psi : \Omega \rightarrow \mathbb{R}$ defined by $\psi(\bm{x}) = \bm{g} \cdot
\bm{x} + c$, with $\bm{g} \in \mathbb{R}^d$ and $c \in \mathbb{R}$, assuming that $u_\sigma =
\psi(\bm{x}_\sigma)$ and $u_K = \psi(\bm{x}_K)$, it results that $u_\sigma - u_K = (\bm{x}_\sigma -
\bm{x}_K)^\intercal \bm{g} = (\bm{x}_\sigma - \bm{x}_K)^\intercal \nabla \psi$. Hence, expression
\eqref{eq:discrete_gradient} leads to $\nabla_K u = \nabla \psi$, which amounts to linear exactness for any
affine function $\psi$ on $K \in \mathcal{M}$, provided that $I_\sigma u = u_\sigma$, which is verified
whenever $\bm{x}_\sigma = \bm{y}_\sigma$, $\forall \sigma \in \mathcal{F}_{int}$, see Figure \ref{fig:meshes}.

Finally, if the face gradient $\nabla_\sigma u$ in Eq.\eqref{eq:face_gradient_corrected} is computed using a
linear interpolation operator applied to the cell gradients reconstructed via the Gauss scheme
\eqref{eq:gauss_discrete_gradient} from both cells sharing the face $\sigma$, the non-orthogonal correction
term in Eq.\eqref{eq:face_gradient_corrected} is associated to a large stencil which includes, besides cells
$K$ and $L$ sharing face $\sigma$, all their neighbouring cells $M \in \mathcal{N}_K \cup \mathcal{N}_L$.

By applying the Gauss corrected scheme from Eq.\eqref{eq:face_gradient_corrected} to the diffusion problem
\eqref{eq:isotropic_diffusion_problem}, one obtains the finite volume scheme
\begin{equation}
  \sum_{L \in \mathcal{N}_K} F_{K,L} + \sum_{\sigma \in \mathcal{F}_{K,ext}} F_{K,\sigma} = \int_K f(\bm{x}) \de\bm{x},
  \quad \forall K \in \mathcal{M}
 \label{eq:fvm_isotropic_tensor_diffusivity_problem}
\end{equation}
where the diffusive fluxes $F_{K,L} = - F^{(d)}_{K,\sigma}$ take the forms
\begin{subequations}
 \begin{align}
   &\begin{aligned}
     F_{K,L} &= \alpha_{K|L} \frac{ | \sigma | }{d_{K,L}} (u_K - u_L)\\
            &- \alpha_{K|L} | \sigma | \left( \bm{n}_{K,\sigma}   
               -\frac{\bm{i}_{K,L}}{\;\bm{n}_{K,\sigma} \cdot \bm{i}_{K,L}} \right) \cdot \nabla_{\sigma} u,
                \;\forall K | L \!\in\! \mathcal{F}_{int}
    \end{aligned}\\
   &\begin{aligned}
     F_{K,\sigma} &= \alpha_{K,\sigma} \frac{ | \sigma | }{d_{K,\sigma}} u_K\\
                &- \alpha_{K,\sigma} | \sigma | \left( \bm{n}_{K,\sigma} 
               - \frac{\bm{i}_{K,\sigma}}{\;\bm{n}_{K,\sigma} \cdot \bm{i}_{K,\sigma}} \right) \cdot \nabla_{\sigma} u, 
                  \;\forall \sigma \!\in\! \mathcal{F}_{K,ext}
    \end{aligned}
 \end{align}
\label{eq:gauss_corrected_fluxes}%
\end{subequations}
with the shorthand notation $d_{K,L} = d_{K,\sigma} + d_{L,\sigma}$ and using the unit vectors
\begin{subequations}
 \begin{align}
   \bm{i}_{K,L} &= \frac{ \bm{x}_L - \bm{x}_K }{ | \bm{x}_L - \bm{x}_K | },
    \quad \forall\, K | L \in \mathcal{F}_{int} \\
   \bm{i}_{K,\sigma} &= \frac{\vphantom{A^A} \bm{x}_\sigma - \bm{x}_K }{ | \bm{x}_\sigma - \bm{x}_K | },
    \quad \forall \sigma \in \mathcal{F}_{K,ext}.
 \end{align}
 \label{eq:cell_to_cell_unit_vector}%
\end{subequations}
The diffusivity in  Eq.\eqref{eq:gauss_corrected_fluxes} is defined as 
\begin{subequations}
 \begin{align}
  \alpha_{K|L} &= \frac{1}{|D_{\sigma}|} \int_{D_\sigma} \alpha(\bm{x}) \,\de \bm{x}, \quad \forall\, K | L \in \mathcal{F}_{int} \\
  \alpha_{K,\sigma} &= \frac{1}{|D_{K,\sigma}|} \int_{D_{K,\sigma}} \alpha(\bm{x}) \,\de \bm{x}, \quad \forall \sigma \in \mathcal{F}_{K,ext},
 \end{align}
 \label{eq:face_diffusivity}%
\end{subequations}
which define piecewise constant functions over the diamond cells $D_\sigma$ and $D_{K,\sigma}$ dual to
internal and external mesh faces respectively. Finally, it is important to notice that the fluxes
\eqref{eq:gauss_corrected_fluxes} are locally conservative, since
\begin{equation}
  F_{K,L} = - F_{L,K}, \quad \forall K|L \in \mathcal{F}_{int}.
 \label{eq:diffusive_flux_conservativity}
\end{equation}

In the definition of the fluxes, a reconstruction of the face gradient $\nabla_{\sigma} u$ must be
employed. For this purpose, a linear interpolation operator is selected at internal faces
\begin{equation}
  \nabla_{\sigma} u = I_\sigma \nabla u
                  = \frac{d_{L,\sigma}}{d_{K,L}} \nabla_K u + \frac{d_{K,\sigma}}{d_{K,L}} \nabla_L u,
   \quad \forall \sigma \in \mathcal{F}_{int} ,
 \label{eq:surface_gradient_interpolation}
\end{equation}
while at boundary faces the simplest choice is $\nabla_{\sigma} u = \nabla_K u$, for all $\sigma \in
\mathcal{F}_{K,ext}$. Here, we will use the approximation
\begin{equation*}
  \nabla_{\sigma} u = -\frac{u_{K}}{d_{K,\sigma}} \bm{n}_{K,\sigma}
                    + \left( \nabla_K u - \left( \bm{n}_{K,\sigma} \cdot \nabla_K u \right) \bm{n}_{K,\sigma} \right),
   \quad \forall \sigma \in \mathcal{F}_{ext} ,
\end{equation*}
in order to recover $\bm{n}_{K,\sigma} \cdot \nabla_{\sigma} u = -u_K/d_{K,\sigma}$ at boundaries.  The linear
interpolation makes use of \eqref{eq:gauss_discrete_gradient} with linear interpolation of the face values
\begin{equation}
  I_\sigma u = \frac{d_{L,\sigma}}{d_{K,L}} u_K + \frac{d_{K,\sigma}}{d_{K,L}} u_L,
   \quad \forall \sigma \in \mathcal{F}_{int}
 \label{eq:linear_interpolated_face_values_in_gauss_gradient}
\end{equation}
while the boundary face values follow directly from the homogeneous Dirichlet conditions in
problem \eqref{eq:isotropic_diffusion_problem}. Additionally, in practical implementations it is customary to
compute the scalar diffusivity $\alpha_{K|L}$ at internal faces $K|L \in \mathcal{F}_{int}$ by a linear
interpolation operator as $\alpha_{K|L} = I_\sigma \alpha$.

To allow for the treatment of non-orthogonal polyhedral meshes, it is useful to consider the associated
isotropic diffusion problem  
\begin{equation}
  \Gamma_\alpha = \alpha \bm{I} ,
 \label{eq:isotropic_diffusivity_tensor}
\end{equation}
where $\Gamma_\alpha$ is isotropic diffusivity tensor associated to the scalar diffusivity $\alpha$. This
allows to reformulate  problem \eqref{eq:isotropic_diffusion_problem}  as
\begin{subequations}
 \begin{align}
  - \div(\Gamma_\alpha \nabla \underline{u}) &= f, \quad \text{in}\; \Omega,\\
  \underline{u} &= 0, \quad \text{on}\; \dd \Omega
 \end{align}
 \label{eq:isotropic_tensor_diffusivity_problem}%
\end{subequations}
with $\Gamma_\alpha(\bm{x})$ naturally verifying the usual assumptions \cite{eymard:2006b}. Similarly, the
associated weak formulation is given by
\begin{equation}
  \begin{aligned}
   &\underline{u} \in H_0^1(\Omega),\\
   &\int_{\Omega} \Gamma_\alpha(\bm{x}) \, \nabla \underline{u}(\bm{x}) \cdot \nabla v(\bm{x}) \,\de\bm{x}
    = \int_{\Omega} f(x) \, v(x) \,\de\bm{x},\\
   &\hspace{18.0em}\quad \forall v \in H_0^1(\Omega).
  \end{aligned}
 \label{eq:isotropic_diffusion_weak_problem}
\end{equation}
It is possible to derive a finite volume scheme for diffusion problems with tensorial diffusivity by
constructing a local discrete gradient \cite{eymard:2006b}, in order to obtain at cell face $\sigma$ a
consistent approximation of the diffusive flux $ -\int_\sigma \left( \Gamma_\alpha(\bm{x}) \nabla
\underline{u}(\bm{x}) \right) \cdot \bm{n}_\sigma \de \gamma(\bm{x})$, with usual notation for finite volume 
schemes.
To this purpose, it is beneficial to rewrite the diffusive flux for an internal face $K|L \in
\mathcal{F}_{int}$ using the diffusivity tensor from Eq.\eqref{eq:isotropic_diffusivity_tensor}. Since
$\Gamma_\alpha$ is symmetric, it follows that
\begin{equation}
  F_{K,L} = -|\sigma| \left( \Gamma_{K,L} \nabla_{K,L} u \right) \cdot \bm{n}_{K,\sigma}
         = -|\sigma| \nabla_{K,L} u  \cdot \left( \Gamma_{K,L} \bm{n}_{K,\sigma} \right).
 \label{eq:isotropic_tensor_symmetry_fvm_fluxes}
\end{equation}
In order to allow for the treatment of non-orthogonal meshes, the following diffusivity tensor decomposition
can be applied
\begin{equation}
  \Gamma_{K,L} = \alpha_{K,L} \bm{I} = \Gamma_{K,L}^{\parallel} + \Gamma_{K,L}^{\nparallel}
 \label{eq:isotropic_diffusivity_decomposition}
\end{equation}
with anisotropic (directional) diffusivity tensors
\begin{subequations}
 \begin{align}
  \Gamma_{K,L}^{\parallel} &= \alpha_{K,L} \, \frac{1}{\left( \bm{i}_{K,L}^\intercal \bm{n}_{K,\sigma} \right)^2}
                         \,\bm{i}_{K,L} \bm{i}_{K,L}^\intercal \\
  \Gamma_{K,L}^{\nparallel} &= \alpha_{K,L} \left( \bm{I}
        - \frac{1}{\left( \bm{i}_{K,L}^\intercal \bm{n}_{K,\sigma} \right)^2} \,\bm{i}_{K,L} \bm{i}_{K,L}^\intercal \right)
 \end{align}
 \label{eq:anisotropic_diffusivity_tensors}%
\end{subequations}
by following the natural directions locally identified from the non-orthogonal polyhedral mesh. Notice also
that both diffusivity tensors are symmetric, since $\Gamma_{K,L}^{\parallel} =
(\Gamma_{K,L}^{\parallel})^\intercal$ and $\Gamma_{K,L}^{\nparallel} = (\Gamma_{K,L}^{\nparallel})^\intercal$,
and   that
\begin{equation}
  \Gamma_{K,L}^{\parallel} = \Gamma_{L,K}^{\parallel} \quad\text{and}\quad \Gamma_{K,L}^{\nparallel} = \Gamma_{L,K}^{\nparallel}.
 \label{eq:anistropic_diffusivity_tensors_property}
\end{equation}
A similar flux decomposition can be carried out at boundary faces $\sigma \in \mathcal{F}_{ext}$ by
substituting the unit vector $\bm{i}_{K,L}$ with $\bm{i}_{K,\sigma}$.

By substituting the tensor decomposition from
Eqs.\eqref{eq:isotropic_diffusivity_decomposition}-\eqref{eq:anisotropic_diffusivity_tensors} into the finite
volume fluxes \eqref{eq:isotropic_tensor_symmetry_fvm_fluxes}, one obtains that
\begin{equation}
  \begin{aligned}
  \Gamma_{K,L} \bm{n}_{K,\sigma} &= \Gamma_{K,L}^{\parallel} \bm{n}_{K,\sigma} + \Gamma_{K,L}^{\nparallel} \bm{n}_{K,\sigma}\\
     &= \alpha_{K,L} \, \frac{1}{\,\bm{i}_{K,L}^\intercal \bm{n}_{K,\sigma}} \,\bm{i}_{K,L}\\
     &\hspace{1.0em}+  \alpha_{K,L} \left( \bm{n}_{K,\sigma} - \frac{1}{\,\bm{i}_{K,L}^\intercal \bm{n}_{K,\sigma}} \,\bm{i}_{K,L} \right),
  \end{aligned}
 \label{eq:tensor_splitting_fvm_fluxes}
\end{equation}
which directly corresponds to the terms of the \emph{Gauss corrected} scheme appearing in
Eq.\eqref{eq:gauss_corrected_fluxes}.  In particular, the first term in
Eq.\eqref{eq:tensor_splitting_fvm_fluxes}, corresponding to the anisotropic diffusivity tensor
$\Gamma_{K,L}^{\parallel}$, is amenable to approximation by a two-point flux scheme, in a manner similar to
what is done in the perpendicular bisection method in \cite{heinemann:1991}. This term, when inserted into the
finite volume diffusive flux, yields
\begin{equation}
  \begin{aligned}
   F_{K,L}^{\parallel} &= -|\sigma| \left( \Gamma_{K,L}^\parallel \nabla_{K,L} u \right) \cdot \bm{n}_{K,\sigma}
         = -|\sigma| \nabla_{K,L} u  \cdot \left( \Gamma_{K,L}^\parallel \bm{n}_{K,\sigma} \right)\\
         &= -\alpha_{K,L} \,\frac{|\sigma|}{\,\bm{i}_{K,L}^\intercal \bm{n}_{K,\sigma}} \bm{i}_{K,L} \cdot \nabla_{K,L} u\\
         &= -\alpha_{K|L} \,\frac{|\sigma|}{\,\bm{i}_{K,L}^\intercal \bm{n}_{K,L}} \frac{u_K-u_L}{|\bm{x}_L-\bm{x}_K|},
  \end{aligned}
 \label{eq:fvm_diffusive_flux_directional_derivative_part}
\end{equation}
that generates a directional derivative which can be easily approximated via a two-point flux scheme. On the
other hand, the second term in the diffusive flux corresponding to the anisotropic diffusivity tensor
$\Gamma_{K,L}^{\nparallel}$ must be treated via a reconstruction of the cell gradient.  it is important also
to notice that, from the tensor decomposition in
Eqs.\eqref{eq:isotropic_diffusivity_decomposition}-\eqref{eq:anisotropic_diffusivity_tensors}, the two-point
flux portion increases its dominance for increasing mesh non-orthogonality, due to the increasing angle
between the unit vectors $\bm{n}_{K,\sigma}$ and $\bm{i}_{K,L}$. This property is beneficial in guaranteeing
diagonal dominance of the linear system matrix and thus numerical stability, as will be clear from the rest of
the discussion.

\section{Discrete weak formulation}
 \label{sec:underlying_weak_form}

Returning to the diffusive fluxes from \emph{Gauss corrected} scheme \eqref{eq:gauss_corrected_fluxes}, using
the diffusivity tensor decomposition in
Eqs.\eqref{eq:isotropic_diffusivity_decomposition}-\eqref{eq:anisotropic_diffusivity_tensors}, the finite
volume fluxes can be rewritten in the form
\begin{subequations}
 \begin{align}
  F_{K,L} &= \alpha_{K|L} \,\tau_{K|L} \, (u_K - u_L)
            - |\sigma| \nabla_{K|L} u \cdot \left( \Gamma_{K,L}^\nparallel \bm{n}_{K,\sigma} \right),\\
  F_{K,\sigma} &= \alpha_{K,\sigma} \,\tau_{K,\sigma} \, u_K
                - |\sigma| \nabla_{K,\sigma} u \cdot \left( \Gamma_{K,\sigma}^\nparallel \bm{n}_{K,\sigma} \right),
 \end{align}
 \label{eq:gauss_corrected_fluxes_weak_form}%
\end{subequations}
for the internal and external faces, respectively. In these formulae, the transmissivities
\begin{equation}
  \tau_{K|L} = \frac{|\sigma|}{d_{K,L}}, \; \forall K | L \in \mathcal{F}_{int}
  \quad \text{and} \quad
  \tau_{K,\sigma} = \frac{|\sigma|}{d_{K,\sigma}}, \; \forall \sigma \in \mathcal{F}_{K,ext}
 \label{eq:transmissibilities}
\end{equation}
are introduced to simplify the notation and $\nabla_{K|L} u$ denotes a generic discrete gradient operator,
still to be defined, that is piecewise constant on the diamond cells $D_{K|L}$ for all $K \in \mathcal{M}$
and $L \in \mathcal{N}_K$.
If one defines the diamond cell gradient from the linear interpolation of cell gradients as in
Eq.\eqref{eq:surface_gradient_interpolation}, then the fluxes become
\begin{subequations}
 \begin{align}
  &\begin{aligned}
   F_{K,L} &= \alpha_{K|L} \,\tau_{K|L} \, (u_K - u_L)\\
          &\hspace{1.0em} + \left( - \nabla_K u \cdot \left( \Gamma_{K,L}^\nparallel \bm{a}_{K,L} \right)
                       + \nabla_L u \cdot \left( \Gamma_{K,L}^\nparallel \bm{a}_{L,K} \right) \right),\\
  \end{aligned}\\
   &F_{K,\sigma} = \alpha_{K,\sigma} \,\tau_{K,\sigma} \, u_K
                - \nabla_K u \cdot \left( \Gamma_{K,\sigma}^\nparallel \bm{a}_{K,\sigma} \right),
 \end{align}
 \label{eq:gauss_corrected_fluxes_weak_form_interpolated}%
\end{subequations}
where the vector quantities
\begin{subequations}
 \begin{align}
  \bm{a}_{K,L} &= |\sigma| \,\frac{d_{L,\sigma}}{d_{K,L}}\, \bm{n}_{K,\sigma},
      \quad  \forall\, K | L \in \mathcal{F}_{int} \\
  \bm{a}_{K,\sigma} &= |\sigma| \,\bm{n}_{K,\sigma},
      \quad  \forall \sigma \in \mathcal{F}_{K,ext}. 
 \end{align}
 \label{eq:a_vector_coeffs}%
\end{subequations}
have been introduced, which are such that $\bm{a}_{K,L} \neq \bm{a}_{L,K}$ generally. Notice also the
approximation introduced in the boundary term $\nabla_{K,\sigma} u \approx \nabla_K u$.  It is now possible to
derive the weak formulation underlying the finite volume scheme
\eqref{eq:fvm_isotropic_tensor_diffusivity_problem}. By multiplying
Eq.\eqref{eq:fvm_isotropic_tensor_diffusivity_problem} by the test function $v_K$ and summing the result for
all $K \in \mathcal{M}$, one obtains
\begin{equation*}
  \sum_{K \in \mathcal{M}} v_K \sum_{L \in \mathcal{N}_K} F_{K,L}
  + \sum_{K \in \mathcal{M}} v_K \sum_{\sigma \in \mathcal{F}_{K,ext}} F_{K,\sigma}
  = \sum_{K \in \mathcal{M}} v_K \int_K f(\bm{x}) \de \bm{x}
\end{equation*}
that, after discrete integration by parts, produces
\begin{multline*}
  \sum_{K|L \in \mathcal{F}_{int}} \left( F_{K,L} \,v_K + F_{L,K} \,v_L \right)
  + \sum_{K \in \mathcal{M}} \sum_{\sigma \in \mathcal{F}_{K,ext}} F_{K,\sigma} v_K\\
  = \sum_{K \in \mathcal{M}} v_K \int_K f(\bm{x}) \de \bm{x}
\end{multline*}
from which, due to flux conservativity \eqref{eq:diffusive_flux_conservativity}, one obtains that
\begin{equation}
 \begin{aligned}
  \sum_{K|L \in \mathcal{F}_{int}} F_{K,L} (v_K - v_L)
  + \sum_{K \in \mathcal{M}} &\sum_{\sigma \in \mathcal{F}_{K,ext}} F_{K,\sigma} v_K\\
  &= \sum_{K \in \mathcal{M}} v_K \int_K f(\bm{x}) \de \bm{x}.
 \end{aligned}
 \label{eq:towards_weak_form_general_case}
\end{equation}

By substituting into Eq.\eqref{eq:towards_weak_form_general_case} the fluxes
\eqref{eq:gauss_corrected_fluxes_weak_form_interpolated} with the face gradient from the linear interpolation
of the  Gauss scheme \eqref{eq:gauss_discrete_gradient}, it is possible to identify two terms $T_1$ and
$T_2 = T_{2,int} + T_{2,ext}$ in the expression
\begin{equation}
  T_1 + T_{2,int} + T_{2,ext} = \sum_{K \in \mathcal{M}} v_K \int_K f(\bm{x}) \de \bm{x},
 \label{eq:towards_weak_form_cell_gradient}
\end{equation}
where
\begin{subequations}
  \begin{align}
   &\hspace{1.25em}\begin{aligned}
    T_1 &= \sum_{K|L \in \mathcal{F}_{int}} \alpha_{K|L} \tau_{K|L} (u_K - u_L)(v_K - v_L)\\
           &\hspace{6.0em}+ \sum_{K \in \mathcal{M}} \sum_{\sigma \in \mathcal{F}_{K,ext}} \alpha_{K,\sigma} \tau_{K,\sigma} u_K v_K,
    \end{aligned}\\
   &\begin{aligned}
    T_{2.int} &= \sum_{K|L \in \mathcal{F}_{int}}\left( -\nabla_K u \cdot \left( \Gamma_{K,L}^\nparallel \bm{a}_{K,L} \right)
                \right.\\
               &\hspace{6.0em}+\left.
                \nabla_L u \cdot \left( \Gamma_{K,L}^{\nparallel} \bm{a}_{L,K} \right) \right) (v_K - v_L),
    \end{aligned}\\
   &\begin{aligned}
    T_{2.ext} = - \sum_{K \in \mathcal{M}} \sum_{\sigma \in \mathcal{F}_{K,ext}} \nabla_K u \cdot
                 \left( \Gamma_{K,\sigma}^\nparallel \bm{a}_{K,\sigma} \right) v_K.
    \end{aligned}
  \end{align}
\end{subequations}
Notice that term $T_1$ defines a symmetric bilinear form
\begin{equation}
 \begin{aligned}
  [u,v]_{\mathcal{D},\alpha,\parallel} &= \sum_{K|L \in \mathcal{F}_{int}} \alpha_{K|L} \tau_{K|L} (u_K - u_L)(v_K - v_L)\\
    &\hspace{6.0em}+ \sum_{K \in \mathcal{M}} \sum_{\sigma \in \mathcal{F}_{K,ext}} \alpha_{K,\sigma} \tau_{K,\sigma} u_K v_K
 \end{aligned}
 \label{eq:inner_product_parallel}
\end{equation}
which is a discretization of the term $\int_\Omega \nabla \underline{u}(\bm{x}) \cdot \left(
\Gamma_\alpha^\parallel(\bm{x}) \nabla v(\bm{x}) \right) \,\de\bm{x}$, directly corresponding to the portion
of diffusive fluxes that can be ascribed to the anisotropic diffusion tensor $\Gamma_\alpha^\parallel$. This
expresses the flux component that is parallel to the local vector $(\bm{x}_L - \bm{x}_K)$ associated to
internal faces $\sigma = K|L$, or to the vector $(\bm{x}_\sigma - \bm{x}_K)$ associated to boundary faces. The
term $T_2 = T_{2,int} + T_{2,ext}$ contains instead the vectors $\bm{a}_{K,L}$ and $\bm{a}_{K,\sigma}$,
related to the diffusivity tensor $\Gamma_\alpha^\nparallel$, but it is not yet in a form readily
corresponding to a discrete weak formulation. To this purpose, it is convenient to rewrite $T_2$ as
\begin{equation*}
 \begin{aligned}
  T_2 &= \sum_{K \in \mathcal{M}} \sum_{L \in \mathcal{N}_K} -\nabla_K u \cdot
                                                   \left(\Gamma_{K,L}^\nparallel \bm{a}_{K,L} \right) (v_K - v_L)\\
    &\hspace{6.0em}- \sum_{K \in \mathcal{M}} \sum_{\sigma \in \mathcal{F}_{K,ext}} \nabla_K u \cdot
                                \left( \Gamma_{K,\sigma}^\nparallel \bm{a}_{K,\sigma} \right) v_K \\
  &= \sum_{K \in \mathcal{M}} \nabla_K u \cdot
                               \left( \sum_{L \in \mathcal{N}_K} \Gamma_{K,L}^\nparallel \bm{a}_{K,L} (v_L - v_K) \right)\\
    &\hspace{6.0em}- \sum_{K \in \mathcal{M}} \nabla_K u \cdot 
                \left(\sum_{\sigma \in \mathcal{F}_{K,ext}} \Gamma_{K,\sigma}^\nparallel \bm{a}_{K,\sigma} v_K \right).
 \end{aligned}
\end{equation*}
The two summation terms between brackets contained in the last expression correspond to the internal faces and
the boundary faces contributions,  respectively. They can be interpreted as a discretization of the term $\int_K
\Gamma_\alpha^\nparallel \nabla v(\bm{x}) \,\de\bm{x}$ for all $K \in \mathcal{M}$, which allows to introduce
the piecewise constant function $(\Gamma_\alpha^\nparallel \nabla v)_{\mathcal{D}}$ that is defined on each
cell $K \in \mathcal{M}$ as
\begin{equation}
 \begin{aligned}
  (\Gamma_\alpha^\nparallel \nabla v)_K
     = \frac{1}{|K|} &\left( \sum_{L \in \mathcal{N}_K} \Gamma_{K,L}^\nparallel \bm{a}_{K,L} (v_L - v_K) \right. \\
                 &\hspace{5.0em}\left.- \sum_{\sigma \in \mathcal{F}_{K,ext}} \Gamma_{K,\sigma}^\nparallel \bm{a}_{K,\sigma} v_K \right),
 \end{aligned}
 \label{eq:inner_product_gamma_npar_nabla_v}
\end{equation}
expressing the fact that the test function gradient cannot be separated from the diffusivity tensor
$\Gamma_\alpha^\nparallel$, since the latter is a face-based quantity defined from the local mesh
non-orthogonality (i.e., from the angle between $\bm{i}_{K,L}$ and $\bm{n}_{K,\sigma}$ unit vectors).  In this
case, the term $T_2$ can be rewritten as a non-symmetric discrete bilinear form
\begin{equation}
  \langle \nabla u ,\nabla v \rangle_{\mathcal{D}.\alpha,\nparallel}
      = \sum_{K\in\mathcal{M}} |K| \nabla_K u \cdot (\Gamma_\alpha^\nparallel \nabla v)_K
      = T_2
 \label{eq:bilinear_form_nparallel}
\end{equation}
which is thus associated to the diffusivity tensor $\Gamma_{\alpha}^{\nparallel}$, expressing the contribution
of diffusion from a local direction not aligned with $(\bm{x}_L - \bm{x}_K)$ at internal faces, or with
$(\bm{x}_\sigma - \bm{x}_K)$ at boundary faces.

Thus, the discrete weak formulation implied when using the \emph{Gauss corrected} scheme for the heterogeneous
isotropic diffusion problem \eqref{eq:isotropic_diffusion_problem} takes the form
\begin{equation}
  \begin{aligned}
   &u \in H_\mathcal{D},\\
   &[u,v]_{\mathcal{D},\alpha,\parallel} + \langle \nabla u ,\nabla v \rangle_{\mathcal{D},\alpha,\nparallel}
    = \sum_{K \in \mathcal{M}} v_K \int_{K} f(\bm{x}) \,\de\bm{x},\\
   &\hspace{20.0em} \forall v \in H_\mathcal{D}.
  \end{aligned}
 \label{eq:discrete_isotropic_diffusion_weak_problem_cell_gradient}
\end{equation}
%

Several remarks are in order on the basis of the previously introduced formulation.
%
  Similarly to \cite{eymard:2006b}, cell gradients can lead to a discrete inner product whenever the mesh
  geometry allows for a direct estimation of face normal fluxes, e.g., in the case of an orthogonal polyhedral
  mesh.
  When instead anisotropic effects (directional bias) emerge locally on cell faces due to mesh
  non-orthogonality, the construction of face gradients becomes inevitable, as done in \cite{eymard:2010}. In
  this latter case, the diffusivity tensor is necessarily defined on diamond cell support.

  Secondly, the linear interpolation operator used to obtain the face gradient $\nabla_\sigma u$ in the
  \emph{Gauss corrected} scheme from a linear combination of cell gradients calculated via the  Gauss
  gradient scheme for all $K|L \in \mathcal{F}_{int}$ implies that
  \begin{equation}
    \nabla_{K,L} u = \frac{|D_{L,\sigma}|}{|D_\sigma|} \nabla_K u + \frac{|D_{K,\sigma}|}{|D_\sigma|} \nabla_L u,
  \end{equation}
  which defines the face gradient as the diamond cell gradient obtained via an inverse volume weighting
  procedure. The same conclusion is also valid for the scalar diffusivity $\alpha_{K|L}$ defined from
  Eq.\eqref{eq:face_diffusivity}.

  Finally, when the scalar diffusivity $\alpha_{K|L}$ appearing inside the fluxes
  \eqref{eq:gauss_corrected_fluxes_weak_form} is computed by a linear interpolation procedure, the
  anisotropic diffusivity tensor $\Gamma_{K,L}^\nparallel$ in the \emph{Gauss corrected} scheme becomes
  \begin{equation*}
   \begin{aligned}
    \Gamma_{K|L}^\nparallel &= \left( \frac{d_{L,\sigma}}{d_{K,L}} \alpha_K + \frac{d_{K,\sigma}}{d_{K,L}} \alpha_L \right)
    \left(\bm{I}-\frac{\bm{i}_{K,L} \bm{i}_{K,L}^\intercal}{\left(\bm{i}_{K,L}^\intercal \bm{n}_{K,\sigma}\right)^2} \right)\\
                &= \left( \frac{d_{L,\sigma}}{d_{K,L}} \alpha_K + \frac{d_{K,\sigma}}{d_{K,L}} \alpha_L \right) \bm{J}_{K,L}
   \end{aligned}
  \end{equation*}
  where the face non-orthogonality symmetric tensor $\bm{J}_{K,L}$ has been defined. As a consequence, the
  diffusive flux \eqref{eq:gauss_corrected_fluxes_weak_form} can be recast into the form
  \begin{equation*}
   \begin{aligned}
    F_{K,L} &= \alpha_{K|L} \,\tau_{K|L}\, (u_K - u_L)
            - |\sigma| \left( \frac{d_{L,\sigma}}{d_{K,L}} \alpha_K + \frac{d_{K,\sigma}}{d_{K,L}} \alpha_L \right)\\
           &\hspace{1.0em} \times \left( \frac{d_{L,\sigma}}{d_{K,L}} \nabla_K u
                  + \frac{d_{K,\sigma}}{d_{K,L}} \nabla_L u \right) \cdot \left( \bm{J}_{K,L} \,\bm{n}_{K,\sigma} \right)
   \end{aligned}
  \end{equation*}
  from which, after introducing the vectors
  \begin{subequations}
    \begin{align}
      \bm{b}_{K,L} &= |\sigma|\left(\frac{d_{L,\sigma}}{d_{K,L}}\right)^2 \bm{J}_{K,L} \bm{n}_{K,\sigma}\\
      \bm{b}_{L,K} &= |\sigma|\left(\frac{d_{K,\sigma}}{d_{K,L}}\right)^2 \bm{J}_{L,K} \bm{n}_{L,\sigma},
    \end{align}%
  \end{subequations}
  one obtains that
  \begin{multline}
     F_{K,L} = \alpha_{K|L} \,\tau_{K|L}\, (u_K - u_L) -\alpha_K \nabla_K u \cdot \bm{b}_{K,L}
              +\alpha_L \nabla_L u \cdot \bm{b}_{L,K}\\
              +|\sigma| \left( \alpha_K \nabla_L u + \alpha_L \nabla_K u \right) \cdot
                \left( \frac{d_{K,\sigma} d_{L,\sigma}}{d_{K,L}^2} \bm{J}_{K,L} \bm{n}_{K,\sigma} \right).
   \label{eq:diffusive_flux_with_special_vectors}
  \end{multline}
 Notice that, if the last term vanishes, the same structure of the anisotropic diffusion fluxes from
  \cite{eymard:2006b} is recovered, similarly to the case of cell based diffusion coefficients, but with
  differently defined $\bm{b}_{K,L}$ and $\bm{b}_{L,K}$ vectors. This implies that a weak formulation similar
  to the one   in \cite{eymard:2006b} can also be obtained in this case. Nevertheless, on generally
  non-orthogonal meshes, the last term in Eq.\eqref{eq:diffusive_flux_with_special_vectors} vanishes only when
  $\alpha_K \nabla_L u = - \alpha_L \nabla_K u$, i.e., only on internal faces where the flux is zero. In all
  the other meaningful cases, the last term in Eq.\eqref{eq:diffusive_flux_with_special_vectors} is non zero
  and it is responsible for the cross terms inside the non-orthogonal correction $\langle \nabla u,
  \nabla v \rangle_{\mathcal{D},\alpha,\nparallel}$ appearing in the weak formulation
  \eqref{eq:discrete_isotropic_diffusion_weak_problem_cell_gradient}.

\section{Convergence analysis}
 \label{sec:convergence_analysis}

The term $T_1$ defined in Eq.\eqref{eq:inner_product_parallel} and appearing in the discrete weak formulation
\eqref{eq:discrete_isotropic_diffusion_weak_problem_cell_gradient} exactly corresponds to the symmetric
bilinear form appearing in \cite{eymard:2006b} for isotropic diffusion operators on polyhedral meshes
satisfying the additional orthogonality condition
\begin{equation}
  (\bm{x}_L - \bm{x}_K) \perp \bm{n}_{K,\sigma}.
 \label{eq:orthogonality_cond}
\end{equation}
However, in the present analysis the same inner product corresponds only to the portion of the discrete
bilinear form containing the contribution to the diffusive flux that is parallel to the local mesh direction,
as identified from the cell-to-cell vector $(\bm{x}_L - \bm{x}_K)$.
Formally, it is possible to define the discrete inner product
\begin{equation}
 \begin{aligned}
  [u,v]_{\mathcal{D},\alpha,\parallel} &= \sum_{K|L \in \mathcal{F}_{int}} \alpha_{K|L} \tau_{K|L} (u_K - u_L)(v_K - v_L)\\
            &\hspace{2.0em}+ \sum_{K \in \mathcal{M}} \sum_{\sigma \in \mathcal{F}_{K,ext}} \alpha_{K,\sigma} \tau_{K,\sigma} u_K v_K
 \end{aligned}
 \label{eq:parallel_inner_product}
\end{equation}
from which the associated norm
\begin{equation}
 \| u \|_\mathcal{D} = ([u,u]_{\mathcal{D},1,\parallel})^{1/2}
 \label{eq:parallel_inner_product}
\end{equation}
directly follows, where we have set  $\alpha = 1$. Such norm verifies the discrete Poincar\'{e} inequality
\begin{equation}
  \|w\|_{L^2(\Omega)} \leq \diam (\Omega) \|w\|_\mathcal{D}, \quad \forall w \in H_\mathcal{D}
 \label{eq:poincare_inequality_discrete}
\end{equation}
as from \cite{eymard:2000}. Furthermore, a relative compactness result in $L^2(\Omega)$ also holds.

\begin{lemma}[\cite{eymard:2006b}, Lemma 2.1]
  Let $\Omega$ be a bounded open connected polyhedral subset of $\mathbb{R}^d$, $d \in \mathbb{N}^\star$ and
 let $(\mathcal{D}_n,u_n)_{n\in\mathbb{N}}$ be a sequence of discretizations such that, for all
 $n \in \mathbb{N}$, $\mathcal{D}_n$ is an admissible finite volume mesh in sense of
 Definition \ref{def:polyhedral_mesh} and $u_n \in H_{\mathcal{D}_n}(\Omega)$. Assume that
 $\lim_{n\rightarrow\infty} h_{\mathcal{D}_n} = 0$ and that there exists a constant $C_1 > 0$ such that
 $\|u\|_{\mathcal{D}_n} \leq C_1$, for all $n \in \mathbb{N}$.  Then there exists a subsequence of
 $(\mathcal{D}_n,u_n)_{n\in\mathbb{N}}$, for simplicity denoted again by $(\mathcal{D}_n,u_n)$, and some
 $\overline{u} \in H_0^1(\Omega)$ such that $u_n$ tends to $\underline{u}$ in $L^2(\Omega)$ as
 $n\rightarrow\infty$, and the inequality
 \begin{equation}
   \int_\Omega |\nabla \underline{u}(\bm{x})|^2 \,\de x \leq \lim_{n\rightarrow\infty} \inf \|u_n\|^2_{\mathcal{D}_n}
 \end{equation}
 holds. Furthermore, for all regular functions $L^\infty(\Omega)$, one has also that
 \begin{equation}
  \begin{aligned}
  \lim_{n\rightarrow\infty} [u_n,P_{\mathcal{D}_n} \varphi]_{\mathcal{D}_n,\alpha,\parallel} &=
   \int_\Omega \Gamma_\alpha^\parallel(\bm{x}) \nabla\underline{u}(\bm{x}) \cdot \nabla \varphi(\bm{x}) \,\de\bm{x},\\
   &\hspace{10.0em} \forall \varphi \in C_c^\infty(\Omega).
  \end{aligned}
 \end{equation}
 with $P_{\mathcal{D}}: C(\Omega) \rightarrow H_{\mathcal{D}}(\Omega)$ the projection operator from Section
 \ref{sec:meshes_spaces}.
 \label{lem:relative_compactness}
\end{lemma}

\begin{proof*}
 The proof is similar to the one reported in \cite{eymard:2006b}, which is obtained for orthogonal meshes,
 even if orthogonality is not strictly required, after substitution of the scalar diffusivity $\alpha$ with
 the diffusivity tensor $\Gamma_\alpha^\parallel$.
\end{proof*}

From the discussion leading to the discrete weak
form \eqref{eq:discrete_isotropic_diffusion_weak_problem_cell_gradient}, it is useful to define a discrete
gradient with anisotropic diffusivity biasing, see also Eq.\eqref{eq:inner_product_gamma_npar_nabla_v}.
\begin{definition}[Discrete gradient with $\Gamma_\alpha^\nparallel$ biasing]
 Let $\Omega$ be a bounded open connected polyhedral subset of $\mathbb{R}^d$, $d \in \mathbb{N}^\star$. Let
 $\mathcal{D}$ be an admissible finite volume discretization in sense of
 Definition \eqref{def:polyhedral_mesh}. The discrete gradient with $\Gamma_\alpha^\nparallel$ anisotropic
 biasing $\nabla_{\mathcal{D},\alpha,\nparallel} : H_\mathcal{D} \rightarrow H_\mathcal{D}^d $ is defined for
 any $u \in H_\mathcal{D}$ as the piecewise constant function
 \begin{equation} 
  \begin{aligned}
   &\nabla_{\mathcal{D},\alpha,\nparallel} u (\bm{x}) = (\Gamma_\alpha^\nparallel \nabla u)_K \\
   &\hspace{0.0em}= \frac{1}{|K|} \left( \sum_{L \in \mathcal{N}_K} \Gamma_{K,L}^\nparallel \bm{a}_{K,L} (u_L - u_K)
       - \sum_{\sigma \in \mathcal{F}_{K,ext}} \Gamma_{K,\sigma}^\nparallel \bm{a}_{K,\sigma} \, u_K \right),\\
   &\hspace{0.0em} \quad \text{for a.e. } \bm{x} \in K, \quad \forall K \in \mathcal{M},
  \end{aligned}
  \label{eq:directional_gradient}
 \end{equation}
 where the discrete anisotropic diffusivity tensor $\Gamma_{K,L}^\nparallel$ (and
 $\Gamma_{K,\sigma}^\nparallel$) is defined in \eqref{eq:anisotropic_diffusivity_tensors} and the vector
 quantities $\bm{a}_{K,L}$ (and $\bm{a}_{K,\sigma}$) are defined
 in \eqref{eq:a_vector_coeffs}.
 \label{def:directional_gradient}
\end{definition}

From the diffusivity tensor decomposition in
Eqs.\eqref{eq:isotropic_diffusivity_decomposition}-\eqref{eq:anisotropic_diffusivity_tensors}, it is possible to
split the $\Gamma_\alpha^\nparallel$-biased discrete gradient into two other discrete gradients.

\begin{definition}[Decomposition of $\Gamma_\alpha^\nparallel$-biased discrete gradient]
 Let $\Omega$ be a bounded open connected polyhedral subset of $\mathbb{R}^d$, $d \in \mathbb{N}^\star$ and
 let $\mathcal{D}$ be an admissible finite volume discretization in sense of
 Definition \eqref{def:polyhedral_mesh}. Let $\nabla_{\mathcal{D},\alpha,\nparallel}$ be the
 $\Gamma_\alpha^\nparallel$-biased discrete gradient, as from Definition \ref{def:directional_gradient}, for
 any $u \in H_\mathcal{D}$. Then the $\Gamma_\alpha^\nparallel$-biased discrete gradient can be decomposed
 into the sum of two other discrete gradients
 \begin{equation}
   \nabla_{\mathcal{D},\alpha,\nparallel} u (\bm{x}) = \nabla_{\mathcal{D},\alpha} u (\bm{x})
                                               - \nabla_{\mathcal{D},\alpha,\parallel} u (\bm{x})
  \label{eq:directional_gradient_decomposition}
 \end{equation}
 where
 \begin{equation} 
  \begin{aligned}
   &\nabla_{\mathcal{D},\alpha} u (\bm{x}) = (\alpha \nabla u)_K \\
   &= \frac{1}{|K|} \left( \sum_{L \in \mathcal{N}_K} \alpha_{K|L} \,\tau_{K|L}\, d_{L,\sigma}
                               \bm{n}_{K,\sigma} (u_L - u_K) \right. \\
  &\left.\hspace{4.0cm} -\sum_{\sigma \in \mathcal{F}_{K,ext}} \alpha_{K,\sigma} \,|\sigma|\, \bm{n}_{K,\sigma} \, u_K \right),\\
   & \quad \text{for a.e. } \bm{x} \in K, \quad \forall K \in \mathcal{M},
  \end{aligned}
  \label{eq:diffusivity_weighted_gradient}
 \end{equation}
 represents a diffusivity weighted discrete gradient, while
 \begin{equation} 
  \begin{aligned}
   &\nabla_{\mathcal{D},\alpha,\parallel} u (\bm{x}) = (\Gamma_{\alpha}^{\parallel} \nabla u)_K \\
   &= \frac{1}{|K|} \left( \sum_{L \in \mathcal{N}_K} \alpha_{K|L} \,\tau_{K|L}\, d_{L,\sigma} 
                                \frac{\bm{i}_{K,L}}{\bm{n}_{K,\sigma} \cdot \bm{i}_{K,L}} (u_L - u_K) \right. \\
&\left.\hspace{4.0cm} -\sum_{\sigma \in \mathcal{F}_{K,ext}} \alpha_{K,\sigma} \,|\sigma|\, 
                                            \frac{\bm{i}_{K,L}}{\bm{n}_{K,\sigma} \cdot \bm{i}_{K,L}} \, u_K \right),\\
   & \quad \text{for a.e. } \bm{x} \in K, \quad \forall K \in \mathcal{M},
  \end{aligned}
  \label{eq:diffusivity_weighted_gradient}
 \end{equation}
 can be interpreted as a $\Gamma_\alpha^\parallel$-biased discrete gradient. 
 \label{def:directional_gradient_decomposition}
\end{definition}

For this finite volume diffusion scheme, the mesh regularity is measured by the factor
\begin{multline}
  \tilde\theta_\mathcal{D}
   = \min \left\{ \min \left\{\; \frac{d_{K,\sigma}}{d_{L,\sigma}}, \frac{h_K}{d_{K,\sigma}}, \bm{n}_{K,L}\cdot\bm{i}_{K,L}
     : \sigma \in \mathcal{F}_{int} \;\right\} \right. , \\
     \min \left\{\; \left. \frac{h_K}{d_{K,\sigma}}, \bm{i}_{K,\sigma}\cdot\bm{n}_{K,\sigma}
     : \sigma \in \mathcal{F}_{ext} \;\right\} \;\right\} \,,
 \label{eq:mesh_regularity_nonparallel}
\end{multline}
which expresses bounds in the empirical measures of mesh regularity that will be presented in Section
\ref{sec:numerical_results}. As a first result, one introduces the bound on the $L^2(\Omega)^d$-norm of the
$\Gamma_\alpha^\nparallel$-biased gradient on any element of $H_\mathcal{D}$.

\begin{lemma}[Bound on $\nabla_{\mathcal{D},\alpha,\nparallel} u$]
 Let $\Omega$ be a bounded open connected polyhedral subset of $\mathbb{R}^d$, $d \in \mathbb{N}^\star$. Let
 $\mathcal{D}$ be an admissible finite volume discretization in sense of Definition \ref{def:polyhedral_mesh}
 and let $0 < \theta \leq \tilde\theta_\mathcal{D}$. Then, the exists $C_1$ depending only on $d$, $\alpha$ and
 $\theta$ such that, for all $u \in H_\mathcal{D}$, one has
 \begin{equation}
   \| \nabla_{\mathcal{D},\alpha,\nparallel} u \|_{L^2(\Omega)^d} \leq C_1 \, \|u\|_\mathcal{D}.
  \label{eq:directional_gradient_bound}
 \end{equation}
 \label{lem:directional_gradient_bound}
\end{lemma}

\begin{proof*}
 Let $u \in H_\mathcal{D}$. Similarly as in \cite{eymard:2006b}, one introduces, for all $K \in \mathcal{M}$,
 $L \in \mathcal{N}_K$ and $\sigma = K|L$ the difference quantities $\delta_{K,\sigma} \bm{x} = (\bm{x}_L
 - \bm{x}_K)$ and $\delta_{K,\sigma} u = (u_L - u_K)$, and for all $\sigma \in \mathcal{F}_{K,\sigma}$ the
 quantities $\delta_{K,\sigma} \bm{x} = (\bm{x}_\sigma - \bm{x}_K)$ and $\delta_{K,\sigma} u = - u_K$. Then,
 the inner product norm in \eqref{eq:parallel_inner_product} leads for a given $K \in \mathcal{M}$ to
 \begin{equation*}
  \begin{aligned}
   \|u\|_{\mathcal{D}}^2 &= [u,u]_{\mathcal{D},1,\parallel}\\
        &= \sum_{K|L \in \mathcal{F}_{int}} \tau_{K|L} (u_L-u_K)^2
          + \sum_{K \in \mathcal{M}} \sum_{\sigma \in \mathcal{F}_{K,ext}} \tau_{K,\sigma} (-u_K)^2\\
        &= \sum_{K \in \mathcal{M}} \frac{1}{2} \sum_{L \in \mathcal{N}_K} \tau_{K|L} (\delta_{K,L} u)^2
          + \sum_{K \in \mathcal{M}} \sum_{\sigma \in \mathcal{F}_{K,ext}} \tau_{K,\sigma} (\delta_{K,\sigma} u)^2.
  \end{aligned}
 \end{equation*}
 Then, Definition \ref{def:directional_gradient} leads to
 \begin{equation*}
   |K| (\Gamma_\alpha^\nparallel \nabla u)_K
   = \sum_{\sigma \in \mathcal{F}_{K}} \alpha_\sigma \,\tau_\sigma\, d_{L,\sigma} \left( \bm{n}_{K,\sigma}
                                - \frac{\delta_{K,\sigma} \bm{x} }{ d_{K,L} } \right) \delta_{K,\sigma} u.
 \end{equation*}
 By using the Cauchy-Schwartz inequality, one obtains that
 \begin{equation*}
  \begin{aligned}
   |K|^2 \, \left| (\Gamma_\alpha^\nparallel \nabla u)_K \right|^2 &\leq
   \sum_{\sigma \in \mathcal{F}_K } \tau_\sigma \alpha_\sigma^2 \left| d_{L,\sigma} \left( \bm{n}_{K,\sigma}
                                - \frac{\delta_{K,\sigma} \bm{x} }{ d_{K,L} } \right) \right|^2\\
          &\hspace{10.0em}\times \sum_{\sigma \in \mathcal{F}_K } \tau_\sigma \left( \delta_{K,\sigma} u \right)^2 ,
  \end{aligned}
 \end{equation*}
 from which, by introducing the upper bound for the scalar diffusivity $C_\alpha \geq \alpha_\sigma^2$, for all
 $\sigma \in \mathcal{F}$, and by noticing that, for $\sigma \in \mathcal{F}_K$, one has
 $\tau_\sigma \leq \frac{|\sigma|}{d_{K,\sigma}}$ and that $\delta_{K,\sigma} \bm{x} \leq (\bm{x}_\sigma
 - \bm{x}_K)$, it follows that
 \begin{equation*}
  \begin{aligned}
   |K|^2 \, \left| (\Gamma_\alpha^\nparallel \nabla u)_K \right|^2
   &\leq C_\alpha \sum_{\sigma \in \mathcal{F}_K } d |D_{K,\sigma}| \left| \frac{d_{L,\sigma}}{d_{K,\sigma}} \bm{n}_{K,\sigma}
                                - \frac{d_{L,\sigma}}{d_{K,\sigma}} \frac{\bm{x}_\sigma - \bm{x}_K}{d_{K,L}} \right|^2\\
                 &\hspace{4.0em}\times \sum_{\sigma \in \mathcal{F}_K } \tau_\sigma \left( \delta_{K,\sigma} u \right)^2 \\
   &\leq C_\alpha d \sum_{\sigma \in \mathcal{F}_K } |D_{K,\sigma}|
                 \left| \frac{d_{L,\sigma}}{d_{K,\sigma}} \bm{n}_{K,\sigma} \right|^2
                 \sum_{\sigma \in \mathcal{F}_K } \tau_\sigma \left( \delta_{K,\sigma} u \right)^2 \\
   &\leq C_\alpha\,\frac{d}{\theta^2}\,|K| \sum_{\sigma \in \mathcal{F}_K } \tau_\sigma \left( \delta_{K,\sigma} u \right)^2.
 \end{aligned}
 \end{equation*}
 Finally, after summing over all $K \in \mathcal{M}$, one obtains that
 \begin{equation*}
  \begin{aligned}
   \sum_{K \in \mathcal{M}} |K| \left| (\Gamma_\alpha^\nparallel \nabla u)_K \right|^2
   &\leq C_\alpha\,\frac{d}{\theta^2} \sum_{K \in \mathcal{M}} \sum_{\sigma \in \mathcal{F}_K } \tau_\sigma
                                                                        \left( \delta_{K,\sigma} u \right)^2\\
   &\leq C_\alpha\,\frac{d}{\theta^2} \sum_{K \in \mathcal{M}}
      \left( \sum_{\sigma \in \mathcal{F}_{K,int}} \tau_\sigma \left( \delta_{K,\sigma} u \right)^2 \right.\\
      &\hspace{8.0em}+\left. 2 \sum_{\sigma \in \mathcal{F}_{K.ext}} \tau_\sigma \left( \delta_{K,\sigma} u \right)^2   \right)\\
   &= 2\, C_\alpha\,\frac{d}{\theta^2} \|u\|_{\mathcal{D}}^2
  \end{aligned}
 \end{equation*}
 from which \eqref{eq:directional_gradient_bound} follows with $C_1 = (1/\theta) \sqrt{2 \,C_\alpha \, d}$.
\end{proof*}


It is now possible to state a weak convergence property for the diffusion weighted discrete gradient.

\begin{lemma}[Weak convergence of $\nabla_{\mathcal{D},\alpha} u$]
 Let $\Omega$ be a bounded open connected polyhedral subset of $\mathbb{R}^d$, $d \in \mathbb{N}^\star$. Let
 $\mathcal{D}$ be an admissible finite volume discretization in sense of
 Definition \eqref{def:polyhedral_mesh} and let $0 < \theta \leq \tilde\theta_\mathcal{D}$. Assume that there
 exists $u \in H_\mathcal{D}$ and a function $\underline{u} \in H_0^1(\Omega)$ such that $u$ tends to
 $\underline{u}$ in $L^2(\Omega)$ as $h_\mathcal{D} \rightarrow 0$, while $\|u\|_\mathcal{D}$ remains
 bounded. Then $\nabla_{\mathcal{D},\alpha} u$ weakly converges to $\alpha \nabla \underline{u}$ in
 $L^2(\Omega)^d$ as $h_\mathcal{D} \rightarrow 0$. Additionally,
 $\nabla_{\mathcal{D},\alpha,\parallel} u$ weakly converges to $\Gamma_\alpha^\parallel \nabla \underline{u}$
 as $h_\mathcal{D} \rightarrow 0$.
 \label{lem:directional_gradient_weak_convergence}
\end{lemma}

\begin{proof*}
 In or
 Let $\varphi \in C_c^\infty(\Omega)$. Assume that $h_{\mathcal{D}}$ is small enough that, for all
 $K \in \mathcal{M}$ and $\bm{x} \in K$, if $\varphi(\bm{x}) \neq 0$ then $\mathcal{F}_{K,ext} = \varnothing$.
 Consider the term $T_1^\mathcal{D}$ defined as
 \begin{equation*}
  \begin{aligned}
   T_1^\mathcal{D} &= \int_\Omega P_\mathcal{D} \varphi(\bm{x}) \nabla_{\mathcal{D},\alpha,\nparallel} u(\bm{x}) \,\de\bm{x}
   = \sum_{K \in\ \mathcal{M}} |K| \,\varphi(\bm{x}_K) \, (\Gamma_\alpha^\nparallel \nabla u)_K \\
   &= \sum_{K|L \in \mathcal{F}_{int}} \left(
       \varphi(\bm{x}_K) \,d_{L,\sigma} \, \alpha_{K|L} \, \tau_{K|L}
                            \vphantom{\frac{(\bm{x}_L - \bm{x}_K)}{d_{K,L}}} \right.\\
        &\hspace{6.0em}\times \left. \left( \bm{n}_{K,\sigma}
                 - \frac{(\bm{x}_L - \bm{x}_K)}{d_{K,L}} \right) (u_L - u_K) \right) \\
      &\qquad +\sum_{K|L \in \mathcal{F}_{int}} \left(
       - \varphi(\bm{x}_L) \,d_{K,\sigma} \, \alpha_{K|L} \, \tau_{K|L} 
                            \vphantom{\frac{(\bm{x}_L - \bm{x}_K)}{d_{K,L}}} \right.\\
       &\hspace{6.0em}\qquad \times \left. \left( \bm{n}_{K,\sigma}
                 - \frac{(\bm{x}_L - \bm{x}_K)}{d_{K,L}} \right) (u_K - u_L) \right) \\
   &= \sum_{K|L \in \mathcal{F}_{int}} \left( \varphi(\bm{x}_K) \,d_{L,\sigma} + \varphi(\bm{x}_L) \,d_{K,\sigma} \right)\\
   &\hspace{6.0em} \times \alpha_{K|L} \, \tau_{K|L}
                          \left( \bm{n}_{K,\sigma} - \frac{(\bm{x}_L - \bm{x}_K)}{d_{K,L}} \right)(u_L - u_K)
  \end{aligned}
 \end{equation*}
 in which $\tau_{K|L} \left( \bm{n}_{K,\sigma} - \frac{(\bm{x}_L - \bm{x}_K)}{d_{K,L}} \right)$ defines a sort
 of non-orthogonal transmissivity, while the first term between brackets can be rewritten as
 \begin{equation*}
  \begin{aligned}
   &\varphi(\bm{x}_K) \,d_{L,\sigma} + \varphi(\bm{x}_L) \,d_{K,\sigma} \\
   &\quad =   \vphantom{\frac{\bm{x}_K + \bm{x}_L}{2}}
        \varphi(\bm{x}_K) (\bm{x}_\sigma - \bm{x}_L) \cdot \bm{n}_{L,\sigma}    
       + \varphi(\bm{x}_L) (\bm{x}_\sigma - \bm{x}_K) \cdot \bm{n}_{K,\sigma} \\
   &\quad = \left( \frac{\varphi(\bm{x}_K) + \varphi(\bm{x}_L)}{2} (\bm{x}_L - \bm{x}_K) \right. \\
    &\hspace{3.0em}+\left. ( \varphi(\bm{x}_K) - \varphi(\bm{x}_L) ) \left( \frac{\bm{x}_K + \bm{x}_L}{2} - \bm{x}_\sigma \right)
      \right) \cdot \bm{n}_{K,\sigma}.
  \end{aligned}
 \end{equation*}
 The term $T_1^\mathcal{D}$ can be decomposed into a sum of two terms $T_1^\mathcal{D} = T_2^\mathcal{D} +
 T_3^\mathcal{D}$, where
 \begin{equation*}
  \begin{aligned}
   T_2^\mathcal{D} &= \sum_{K|L \in \mathcal{F}_{int}} \bm{n}_{K,\sigma} \cdot (\bm{x}_L - \bm{x}_K)
                    \frac{\varphi(\bm{x}_K) + \varphi(\bm{x}_L)}{2} \\
                    &\hspace{2.0cm} \times \alpha_{K|L}\, \tau_{K|L}
                    \left( \bm{n}_{K,\sigma} - \frac{(\bm{x}_L - \bm{x}_K)}{d_{K,L}} \right)(u_L - u_K), \\
   T_3^\mathcal{D} &= \sum_{K|L \in \mathcal{F}_{int}} \bm{n}_{K,\sigma} \cdot
                    \left( \frac{\bm{x}_K + \bm{x}_L}{2} - \bm{x}_\sigma \right)
                    ( \varphi(\bm{x}_K) - \varphi(\bm{x}_L) ) \\
                    &\hspace{2.0cm} \times \alpha_{K|L}\, \tau_{K|L}
                    \left( \bm{n}_{K,\sigma} - \frac{(\bm{x}_L - \bm{x}_K)}{d_{K,L}} \right)(u_L - u_K).
  \end{aligned}
 \end{equation*}
 Starting with the analysis of term $T_3^\mathcal{D}$, by Cauchy-Schwartz inequality one gets
 \begin{equation*}
  \begin{aligned}
   \left( T_3^\mathcal{D} \right)^2 \leq
     & \sum_{K|L \in \mathcal{F}_{int}} \tau_{K|L} \, \alpha_{K|L}^2
              \left| \frac{\bm{x}_K + \bm{x}_L}{2} - \bm{x}_\sigma \right|^2\\
     & \qquad \times (\varphi(\bm{x}_K) - \varphi(\bm{x}_L))^2
              \left| \bm{n}_{K,\sigma} - \frac{(\bm{x}_L - \bm{x}_K)}{d_{K,L}} \right|^2\\
     & \qquad \qquad \times \sum_{K|L \in \mathcal{F}_{int}} \tau_{K|L} (u_L - u_K)^2,
  \end{aligned}
 \end{equation*}
 in which, due to triangle inequality
 \begin{equation*}
  \left| \frac{\bm{x}_K + \bm{x}_L}{2} - \bm{x}_\sigma \right|
  \leq \frac{1}{2} |\bm{x}_L - \bm{x}_\sigma| + \frac{1}{2} |\bm{x}_K - \bm{x}_\sigma|
  \leq h_{\mathcal{D}},
 \end{equation*}
 while due to mesh regularity
 \begin{multline*}
  \left| \bm{n}_{K,\sigma} - \frac{(\bm{x}_L - \bm{x}_K)}{d_{K,L}} \right|\\
  = \left| \bm{n}_{K,\sigma} - \frac{\bm{i}_{K,L}}{\bm{n}_{K,\sigma} \cdot \bm{i}_{K,L}} \right|
  \leq 1 + \left|\frac{\bm{i}_{K,L}}{\bm{n}_{K,\sigma} \cdot \bm{i}_{K,L}} \right|
  \leq 1 + \frac{1}{\theta},
 \end{multline*}
 from which, after introducing $C_\alpha \geq \alpha_{K|L}^2$, it follows that
 \begin{equation*}
  \left( T_3^\mathcal{D} \right)^2 \leq
  C_\sigma C_2 \left( 1 + \frac{1}{\theta} \right)^2 h_\mathcal{D}^2 \,|\Omega|\, \|u\|_{\mathcal{D}}^2
 \end{equation*}
 with $C_2$ only depending on $d$, $\Omega$ and $\varphi$. Thus one concludes that
 $\lim_{h_\mathcal{D} \rightarrow 0} T_3^\mathcal{D} = 0$.
 Successively, consider the term $T_2^\mathcal{D}$ that can be rewritten as the sum of two terms
 \begin{equation*}
  \begin{aligned}
   T_2^\mathcal{D} &= \sum_{K|L \in \mathcal{F}_{int}}\alpha_{K|L}\,|\sigma| \frac{\varphi(\bm{x}_K) + \varphi(\bm{x}_L)}{2}\\
         &\hspace{6.0em} \times\left( \bm{n}_{K,\sigma} - \frac{\bm{i}_{K,L}}{\bm{n}_{K,\sigma} \cdot \bm{i}_{K,L}} \right) (u_L - u_K)\\
   &= \sum_{K|L \in \mathcal{F}_{int}} \alpha_{K|L} \, |\sigma| \, \bm{n}_{K,\sigma}
                                            \frac{\varphi(\bm{x}_K) + \varphi(\bm{x}_L)}{2} (u_L - u_K)\\
   &\hspace{2.0em} - \sum_{K|L \in \mathcal{F}_{int}} \alpha_{K|L} \, |\sigma| \,
                                            \frac{\bm{i}_{K,L}}{\bm{n}_{K,\sigma} \cdot \bm{i}_{K,L}}
                                            \frac{\varphi(\bm{x}_K) + \varphi(\bm{x}_L)}{2} (u_L - u_K)\\
   &= T_{2,1}^\mathcal{D} + T_{2,2}^\mathcal{D}.
  \end{aligned}
 \end{equation*}
 Compare term $T_{2,1}^\mathcal{D}$ with the term
 \begin{equation*}
  \begin{aligned}
   T_4^\mathcal{D} &= -\int_\Omega \alpha(\bm{x}) \, u(\bm{x}) \nabla \varphi (\bm{x}) \,\de\bm{x} \\
   &= -\sum_{K \in \mathcal{M}} \sum_{\sigma \in \mathcal{F}_{K,int}} \alpha_{K|L}\, u_K \int_{K|L} \varphi(\bm{x})
                                                               \, \bm{n}_{K,\sigma} \,\de\gamma(\bm{x}) \\
   &= \sum_{K|L \in \mathcal{F}_{int}} \alpha_{K|L} (u_L - u_K) \int_{K|L} \varphi(\bm{x})
                                                               \, \bm{n}_{K,\sigma} \,\de\gamma(\bm{x}),
  \end{aligned}
 \end{equation*}
 which is such that
 \begin{equation*}
  \lim_{h_\mathcal{D} \rightarrow 0} T_4^\mathcal{D}
  = - \int_\Omega \alpha(\bm{x}) \underline{u}(\bm{x}) \nabla\varphi(\bm{x}) \,\de\bm{x}
  =  \int_\Omega \alpha(\bm{x}) \varphi(\bm{x}) \nabla \underline{u}(\bm{x}) \,\de\bm{x}.
 \end{equation*}
 Due to the fact that midpoint face interpolation is first order accurate
 \begin{equation*}
   \left| \frac{1}{|\sigma|} \int_{K|L} \varphi(\bm{x}) \,\de\gamma(\bm{x})
          - \frac{\varphi(\bm{x}_K) + \varphi(\bm{x}_L)}{2} \right|
   \leq h_{\mathcal{D}} \| \nabla \varphi \|_{L^\infty(\Omega)},
 \end{equation*}
 one has that
 \begin{equation*}
  \begin{aligned}
   &\left( T_4^\mathcal{D} - T_{2,1}^\mathcal{D} \right)^2\\
   &\hspace{1.0em}\leq \sum_{K|L \in \mathcal{F}_{int}} (\alpha_{K|L} \,|\sigma|\, \bm{n}_{K,\sigma})^2 (u_L - u_K)^2\\
   &\hspace{3.0em} \times \sum_{K|L \in \mathcal{F}_{int}} \left| \frac{1}{|\sigma|}
        \int_{K|L} \varphi(\bm{x}) \,\de\gamma(\bm{x}) - \frac{\varphi(\bm{x}_K) + \varphi(\bm{x}_L)}{2} \right|^2\\
   &\hspace{1.0em}\leq \sum_{K|L \in \mathcal{F}_{int}} (\alpha_{K|L} \,|\sigma|\, \bm{n}_{K,\sigma})^2 (u_L - u_K)^2
        \sum_{K|L \in \mathcal{F}_{int}}  h_{\mathcal{D}} \| \nabla \varphi \|_{L^\infty(\Omega)},
  \end{aligned}
 \end{equation*}
 from which it follows that $\lim_{h_\mathcal{D} \rightarrow 0} \left( T_4^\mathcal{D} -
 T_{2,1}^\mathcal{D} \right)^2 = 0$.  Thus, $T_{2,1}^\mathcal{D} = T_1^\mathcal{D} - T_{2,2}^\mathcal{D}$
 converges to $T_4^\mathcal{D}$ and, due to density of $C_c^\infty(\Omega)$ in $L^2(\Omega)$,
 $\nabla_{\mathcal{D},\alpha} u$ weakly converges to $\alpha \nabla \underline{u}$ as
 $h_{\mathcal{D}} \rightarrow 0$.
 Thus, the term $T_{2,2}^\mathcal{D}$ can be compared to
 \begin{equation*}
  \begin{aligned}
   T_5^\mathcal{D} &= -\int_\Omega \left( \Gamma_\alpha^\nparallel(\bm{x}) - \alpha(\bm{x}) \bm{I} \right)
                               u(\bm{x}) \nabla \varphi (\bm{x}) \,\de\bm{x} \\
   &= -\sum_{K \in \mathcal{M}} \sum_{\sigma \in \mathcal{F}_{K,\int}}
                 \alpha_{K|L}\, \frac{\bm{i}_{K,L} \bm{i}_{K,L}^\intercal}{\left( \bm{i}_{K,L}^\intercal \bm{n}_{K,\sigma} \right)^2}
                  u_K \int_{K|L} \varphi(\bm{x}) \, \bm{n}_{K,\sigma} \,\de\gamma(\bm{x})
  \end{aligned}
 \end{equation*}
 which is such that
 \begin{equation*}
  \lim_{h_\mathcal{D} \rightarrow 0} T_5^\mathcal{D}
  = \int_\Omega \Gamma_\alpha^\parallel(\bm{x}) \underline{u}(\bm{x}) \nabla\varphi(\bm{x}) \,\de\bm{x}
  = -\int_\Omega \Gamma_\alpha^\parallel(\bm{x}) \varphi(\bm{x}) \nabla \underline{u}(\bm{x}) \,\de\bm{x}.
 \end{equation*}
 By a similar procedure, the term $T_{2,2}^\mathcal{D}$ converges to term $T_5^\mathcal{D}$ and so
 $\nabla_{\mathcal{D},\alpha,\parallel} u$ weakly converges to $\Gamma_{\alpha}^\parallel \nabla
 \underline{u}$ as $h_{\mathcal{D}} \rightarrow 0$.
\end{proof*}

The diffusion weighted discrete gradient provides indeed a consistent gradient scheme.

\begin{lemma}[Consistency of $\nabla_{\mathcal{D},\alpha}$]
 Let $\Omega$ be a bounded open connected polyhedral subset of $\mathbb{R}^d$, $d \in \mathbb{N}^\star$. Let
 $\mathcal{D}$ be an admissible finite volume discretization in sense of
 Definition \eqref{def:polyhedral_mesh} and let $0 < \theta \leq \tilde\theta_\mathcal{D}$. Let
 $\underline{u} \in C^2(\overline{\Omega})$ be such that $\underline{u} = 0$ on $\dd \Omega$. Then there exists
 $C_3$, depending only on $\Omega$, $\theta$, $\underline{u}$ and $\alpha$ such that
 \begin{equation}
   \| \nabla_{\mathcal{D},\alpha} P_{\mathcal{D}} \underline{u} - \alpha \nabla \underline{u} \|_{L^2(\Omega)^d} \leq C_3 h_{\mathcal{D}}
  \label{eq:directional_gradient_consistency}
 \end{equation}
 \label{lem:directional_gradient_consistency}
\end{lemma}

\begin{proof*}
 From Definition \ref{def:directional_gradient_decomposition} for any $K \in \mathcal{M}$ one has
 \begin{equation*}
  \begin{aligned}
    |K| (\nabla_{\mathcal{D},\alpha} P_\mathcal{D} u)_K
   = &\sum_{L \in \mathcal{N}_K} \alpha_{K|L} \,\tau_{K|L}\, d_{L,\sigma} \bm{n}_{K,\sigma}
           (\underline{u}(\bm{x}_L)-\underline{u}(\bm{x}_K)) \\
     &\hspace{2.0cm} -\sum_{K \in \mathcal{M}} \sum_{\sigma \in \mathcal{F}_{K,ext}} \alpha_{K,\sigma} \,|\sigma|\, \bm{n}_{K,\sigma} \underline{u}(\bm{x}_K)
  \end{aligned}
 \end{equation*}
 Let $(\alpha \nabla \underline{u})_K$ be the mean value of $\alpha \nabla \underline{u}$ over $K$
 \begin{equation*}
   (\alpha \nabla \underline{u})_K = \frac{1}{|K|} \int_K \alpha(\bm{x}) \nabla \underline{u}(\bm{x}) \,\de\bm{x}.
 \end{equation*}
 Due to the regularity of $\underline{u}$ and the homogeneous Dirichlet boundary conditions, the flux
 consistency error estimates include a constant $C_4$, only depending on $L^\infty$ norm of second derivatives
 of $\underline{u}$ (and of $\alpha$), such that for all $\sigma = K|L \in \mathcal{F}_{int}$, one has
 \begin{equation*}
   |e_\sigma| \leq C_4 h_{\mathcal{D}} \quad
    \text{with} \quad e_\sigma = (\alpha\nabla \underline{u})_K \cdot {n}_{K,\sigma}
                              - \alpha_{K|L} \frac{\underline{u}(\bm{x}_L) - \underline{u}(\bm{x}_K)}{d_{K,L}}
 \end{equation*}
 while for all $\sigma \in \mathcal{F}_{ext}$ one has
 \begin{equation*}
   |e_\sigma| \leq C_4 h_{\mathcal{D}} \quad
    \text{with} \quad e_\sigma = (\alpha\nabla \underline{u})_K \cdot {n}_{K,\sigma}
                              - \alpha_{K,\sigma} \frac{- \underline{u}(\bm{x}_K)}{d_{K,\sigma}}
 \end{equation*}
 These flux consistency errors allow to recast $(\nabla_{\mathcal{D},\alpha} P_\mathcal{D} u)_K$ as
 \begin{equation*}
  \begin{aligned}
   |K| (\nabla_{\mathcal{D},\alpha} P_\mathcal{D} u)_K
    = &\sum_{L \in \mathcal{N}_K} |\sigma| d_{L,\sigma} \bm{n}_{K,\sigma} (\alpha \nabla \underline{u})_K \cdot \bm{n}_{K,\sigma}\\
      &\hspace{2.0em}
     -\sum_{\sigma \in \mathcal{F}_{K,ext}} |\sigma| d_{K,\sigma} \bm{n}_{K,\sigma} (\alpha \nabla \underline{u})_K \cdot \bm{n}_{K,\sigma}
     + R_K 
  \end{aligned}
 \end{equation*}
 where the consistency residual term is defined as
 \begin{equation*}
   R_K = - \sum_{L \in \mathcal{N}_K} |\sigma| d_{L,\sigma} \bm{n}_{K,\sigma} e_\sigma
         - \sum_{\sigma \in \mathcal{F}_{K,ext}} |\sigma| d_{K,\sigma} \bm{n}_{K,\sigma} e_\sigma.
 \end{equation*}
From the geometrical identity valid for any vector $\bm{x}_0, \bm{v} \in \mathbb{R}^d$ and for all $K \in
 \mathcal{M}$
 \begin{equation}
   \frac{1}{|K|} \sum_{\sigma \in \mathcal{F}_K} |\sigma| (\bm{x}_\sigma - \bm{x}_0) \bm{n}_{K,\sigma} \cdot \bm{v}
   = \bm{v} ,
 \end{equation}
 which is a direct consequence of Eq.\eqref{eq:cell_volume_tensor_identity}, it follows that
 \begin{equation*}
  |K| (\nabla_{\mathcal{D},\alpha} P_\mathcal{D} u)_K \leq \frac{1}{\theta} |K| (\alpha \nabla \underline{u})_K + R_K.
 \end{equation*}
 Due to flux consistency error estimates, it also follows that
 \begin{equation}
  \begin{aligned}
   |R_K| &\leq \sum_{L \in \mathcal{N}_K} |\sigma| d_{L,\sigma} |e_\sigma|
             + \sum_{\sigma \in \mathcal{F}_{K,ext}} |\sigma| d_{K,\sigma} |e_\sigma| \\
         &\leq \frac{C_4}{\theta} h_{\mathcal{D}} \sum_{\sigma \in \mathcal{F}_{K}} |\sigma| d_{K,\sigma}
         = C_4 \frac{d}{\theta}\,|K|\, h_{\mathcal{D}}.
  \end{aligned}
 \end{equation}
 As a consequence, one obtains that
 \begin{equation}
  \begin{aligned}
   &\sum_{K \in \mathcal{M}} |K| \left| (\nabla_{\mathcal{D},\alpha} P_{\mathcal{D}} u)_K - (\alpha \nabla \underline{u})_K \right|^2\\
   &\hspace{2.0em}\leq \sum_{K \in \mathcal{M}} C_4^2 \left( \frac{d}{\theta} \right)^2 h_{\mathcal{D}}^2 |K|
    = \left( C_4 \frac{d}{\theta} \right)^2 h_{\mathcal{D}}^2 |\Omega|.
  \end{aligned}
  \label{eq:density_weighted_gradient_consistency_a}
 \end{equation}
 Using $\underline{u} \in C^2(\Omega)$ and $\alpha$ regularity, there exists $C_5$, only dependent on
 $L^\infty$ norm of the second derivatives of $\underline{u}$, such that
 \begin{equation}
  \sum_{K \in \mathcal{M}} \int_K \left| \alpha \nabla \underline{u} - (\alpha \nabla \underline{u})_K \right|^2
   \leq C_5 h_\mathcal{D}^2.
  \label{eq:density_weighted_gradient_consistency_b}
 \end{equation}
 From Eqs.\eqref{eq:density_weighted_gradient_consistency_a}
 and \eqref{eq:density_weighted_gradient_consistency_b}, one gets the existence of $C_c$, only dependent on
 $\Omega$, $\theta$, $\underline{u}$ and $\alpha$, such that \eqref{eq:directional_gradient_consistency} holds.
\end{proof*}

 In order to complete the convergence analysis, an upper bound in
$\|\cdot\|_{\mathcal{D}}$ and the properties of weak consistency and convergence must also be proved for the
Gauss gradient scheme appearing in the discrete bilinear form defined in
Eq.\eqref{eq:bilinear_form_nparallel}. Remember that, as also argued in \cite{eymard:2009b}, even though
consistent, the Gauss gradient scheme does not allow to obtain coercivity and hence uniqueness. Nevertheless,
in practical implementations also the Gauss gradient scheme suffices in obtaining stable coercive
diffusion operators, as widely verified in the numerical tests reported in Section
\ref{sec:numerical_results}. This is both a consequence of the limited non-orthogonality encountered in
unstructured meshes used herein, and so a limited importance of the non-orthogonal correction term in the
\emph{Gauss corrected} scheme, and also a result of the action of $[u,v]_{\mathcal{D},\alpha}$ term,  that
provides a consistent stabilization term to the bilinear form $\langle u,v
\rangle_{\mathcal{D},\alpha,\nparallel}$.

It is now possible to prove convergence of the weak formulation associated to the \emph{Gauss corrected}
scheme, provided an assumption of the limited contribution from the mesh non-orthogonality correction term.
In particular, one requires the condition of small gradient distortion
\begin{equation}
 \sum_{K \in \mathcal{M}} |K| \nabla_{\mathcal{D}} u \cdot \nabla_{\mathcal{D},1} u
 \geq \sum_{K \in \mathcal{M}} |K| \nabla_{\mathcal{D}} u \cdot \nabla_{\mathcal{D},1,\parallel} u
 \label{eq:gradient_alignment_sufficient_condition_convergence}
\end{equation}
which states that the diffusivity weighted discrete gradient with unit diffusivity $\nabla_{\mathcal{D},1} u$
has a preferential alignment with the  Gauss (stabilized) gradient $\nabla_{\mathcal{D}} u$ with respect
to the $\Gamma_\alpha^\parallel$-weighted gradient $\nabla_{\mathcal{D},1,\parallel} u$ with unit diffusivity
$\alpha=1$. Notice that this sufficient condition for convergence is already known to finite volume
practitioners, which usually require limited mesh non-orthogonality to have stable discretizations. In the
present analysis, the role of condition \eqref{eq:gradient_alignment_sufficient_condition_convergence} is made
clear in proving the discrete $H^1(\Omega)$ estimate.

\begin{lemma}[Discrete $H^1(\Omega)$ estimate]
 Under assumption \eqref{eq:gradient_alignment_sufficient_condition_convergence} and the hypotheses of the
 heterogeneous diffusion problem \eqref{eq:isotropic_diffusion_problem}, let $\mathcal{D}$ be an admissible
 finite volume discretization in sense of Definition \eqref{def:polyhedral_mesh} and let $0
 < \theta \leq \tilde\theta_\mathcal{D}$. Assume that $0 < \alpha_0 \leq \alpha(\bm{x})$ for
 a.e. $\bm{x} \in \Omega$ and also assume that $u \in H_{\mathcal{D}}$ is a solution of the discrete weak
 problem \eqref{eq:discrete_isotropic_diffusion_weak_problem_cell_gradient}. Then the following estimate holds
 \begin{equation}
   \alpha_0 \|u\|_{\mathcal{D}} \leq \diam (\Omega) \|f\|_{L^2(\Omega)}.
  \label{eq:h1_estimate}
 \end{equation}
 \label{lem:h1_estimate}
\end{lemma}

\begin{proof*}
 Consider the discrete weak formulation \eqref{eq:discrete_isotropic_diffusion_weak_problem_cell_gradient}
 and set $v = u,$ to obtain
 \begin{equation*}
   [u,u]_{\mathcal{D},\alpha,\parallel} + \langle \nabla u ,\nabla u \rangle_{\mathcal{D}.\alpha,\nparallel}
   = \int_{\Omega} u(\bm{x}) f(\bm{x}) \,\de\bm{x}.
 \end{equation*}
 After assumption \eqref{eq:gradient_alignment_sufficient_condition_convergence} it follows that
 \begin{equation*}
   \sum_{K \in \mathcal{M}} |K| \nabla_{\mathcal{D}} u \cdot \nabla_{\mathcal{D},1,\nparallel} u \geq 0
 \end{equation*}
 which allows to conclude that
 \begin{equation*}
  \langle \nabla u ,\nabla u \rangle_{\mathcal{D}.\alpha,\nparallel}
  = \sum_{K \in \mathcal{M}} \nabla_K u \cdot (\Gamma_\alpha^\nparallel \nabla u)_K
  \geq \sum_{K \in \mathcal{M}} \alpha_0 \nabla_K u \cdot \nabla_{\mathcal{D},1,\nparallel} u \geq 0 ,
 \end{equation*}
 from which one obtains
 \begin{equation}
   [u,u]_{\mathcal{D},\alpha,\parallel} + \langle \nabla u ,\nabla u \rangle_{\mathcal{D}.\alpha,\nparallel}
    \geq [u,u]_{\mathcal{D},\alpha,\parallel} \geq \alpha_0 \|u\|_{\mathcal{D}}^2.
  \label{eq:discrete_h1_inequality_a}
 \end{equation}
 From the Cauchy-Schwartz inequality and from the discrete Poincar\'{e}
 inequality \eqref{eq:poincare_inequality_discrete} one also obtains
 \begin{equation}
  \begin{aligned}
   \int_{\Omega} u(\bm{x}) f(\bm{x}) \,\de\bm{x}
    &\leq \|u\|_{L^2(\Omega)} \|f\|_{L^2(\Omega)} \\
    &\leq  \diam (\Omega) \|u\|_{\mathcal{D}} \|f\|_{L^2(\Omega)}.
  \end{aligned}
  \label{eq:discrete_h1_inequality_b}
 \end{equation}
 Combining together Eqs.\eqref{eq:discrete_h1_inequality_a} and \eqref{eq:discrete_h1_inequality_b} allows to
 recover the discrete estimate \eqref{eq:h1_estimate}.
\end{proof*}

\begin{corollary}[Existence and uniqueness of a discrete solution]
Assume \eqref{eq:gradient_alignment_sufficient_condition_convergence} and the hypotheses of
 problem \eqref{eq:isotropic_diffusion_problem}. Let $\mathcal{D}$ be an admissible finite volume
 discretization in sense of Definition \eqref{def:polyhedral_mesh} and let
 $0 < \theta \leq \tilde\theta_\mathcal{D}$. Then, there exists a unique solution to
 problem \eqref{eq:discrete_isotropic_diffusion_weak_problem_cell_gradient}.
 \label{cor:existence_uniqueness}
\end{corollary}

\begin{proof*}
 Assume $f=0$ in the finite dimensional
 system \eqref{eq:discrete_isotropic_diffusion_weak_problem_cell_gradient}. From the discrete Poincar\'{e}
 inequality \eqref{eq:poincare_inequality_discrete} one gets $u=0$, thus proving that the linear
 problem \eqref{eq:discrete_isotropic_diffusion_weak_problem_cell_gradient} is uniquely solvable.
\end{proof*}

Finally, it is possible to state the convergence of the finite volume \emph{Gauss corrected} scheme to the solution of the
associated weak problem \eqref{eq:discrete_isotropic_diffusion_weak_problem_cell_gradient}.

\begin{theorem}[Convergence of \emph{Gauss corrected} scheme]
 Assuming the sufficient condition \eqref{eq:gradient_alignment_sufficient_condition_convergence} and the
 hypotheses of the heterogeneous diffusion problem \eqref{eq:isotropic_diffusion_problem}, let $\mathcal{D}$
 be an admissible finite volume discretization in sense of Definition \eqref{def:polyhedral_mesh} with
 $0 < \theta \leq \tilde\theta_\mathcal{D}$. Assume that $0 < \alpha_0 \leq \alpha(\bm{x})$ for
 a.e. $\bm{x} \in \Omega$ and also assume that $u \in H_{\mathcal{D}}$ is a solution to the discrete weak
 problem \eqref{eq:discrete_isotropic_diffusion_weak_problem_cell_gradient}. Then $u$ converges in
 $L^2(\Omega)$ to $\underline{u}$, that is the weak solution to problem \eqref{eq:isotropic_diffusion_problem}
 in the sense of \eqref{eq:isotropic_diffusion_weak_problem}, as $h_{\mathcal{D}} \rightarrow 0$.
 \label{th:convergence_gauss_corrected}
\end{theorem}

\begin{proof*}
 The convergence proof uses the compactness technique presented in \cite{eymard:2007d}.
 Consider a subsequence of admissible discretizations $(\mathcal{D}_n)_{n \in \mathbb{N}}$ such that
 $h_{\mathcal{D}} \rightarrow 0$ as $n \rightarrow \infty$ while $\theta_{\mathcal{D}_n} \geq \theta$ for all
 $n \in \mathbb{N}$. Using Lemma \ref{lem:h1_estimate} one can apply Lemma \ref{lem:relative_compactness},
 which is a discrete counterpart of the Rellich theorem and gives the existence of a subsequence, for
 simplicity denoted again with $(\mathcal{D}_n)_{n \in \mathbb{N}}$, and of some $\underline{u} \in
 H_0^1(\Omega)$ such that the solution $u_{\mathcal{D}_n}$ to
 problem \eqref{eq:discrete_isotropic_diffusion_weak_problem_cell_gradient} tends to $\underline{u}$ in
 $L^2(\Omega)$ as $n \rightarrow \infty$. Let $\varphi \in C_c^\infty(\Omega)$ and select $v =
 P_{\mathcal{D}} \varphi$ as test function in
 problem \eqref{eq:discrete_isotropic_diffusion_weak_problem_cell_gradient}, from which one has
 \begin{equation}
   [u, P_{\mathcal{D}_n} \varphi]_{\mathcal{D}_n,\alpha,\parallel}
    + \langle \nabla u ,\nabla P_{\mathcal{D}_n} \varphi \rangle_{\mathcal{D}_n.\alpha,\nparallel}
    = \int_{\Omega} f(\bm{x}) P_{\mathcal{D}_n} \varphi \,\de\bm{x}.
  \label{eq:convergence_a}
 \end{equation}


 Let then $n \rightarrow \infty$ in Eq.\eqref{eq:convergence_a}. Thanks to Lemma
 \ref{lem:directional_gradient_weak_convergence} and Lemma \ref{lem:directional_gradient_consistency}, but
 also to their counterparts for the $\nabla_\mathcal{D} u$ in  Appendix \ref{sec:appendix}, considering the decomposition
 \begin{equation*}
  \langle \nabla u ,\nabla u \rangle_{\mathcal{D},\alpha,\nparallel}
   = \langle \nabla u ,\nabla u \rangle_{\mathcal{D},\alpha}
     - \langle \nabla u ,\nabla u \rangle_{\mathcal{D},\alpha,\parallel} ,
 \end{equation*}
 one obtains the convergence of the diffusivity weighted gradient portion $\langle \nabla u ,\nabla u
 \rangle_{\mathcal{D},\alpha}$ contained within $\langle \nabla u ,\nabla u
 \rangle_{\mathcal{D}.\alpha,\nparallel}$ term
 \begin{equation*}
   \lim_{n \rightarrow \infty} \int_{\Omega} \nabla_{\mathcal{D}_n} u \cdot (\Gamma_\alpha \nabla P_{\mathcal{D}_n} \varphi) \,\de\bm{x}
    = \int_\Omega \nabla \underline{u} \cdot (\Gamma_\alpha) \nabla \varphi \,\de\bm{x}.
 \end{equation*}
 From Lemma \ref{lem:relative_compactness} and
 Lemma \ref{lem:directional_gradient_weak_convergence} one also
 obtains that the sum of the remaining terms in
 problem \eqref{eq:discrete_isotropic_diffusion_weak_problem_cell_gradient} is such that
 \begin{equation*}
  \lim_{n \rightarrow \infty} \int_{\Omega} \left( [u,P_{\mathcal{D}_n} \varphi]_{\mathcal{D}_n,\alpha,\parallel}
    - \nabla_{\mathcal{D}_n} u \cdot (\Gamma_\alpha^\parallel \nabla P_{\mathcal{D}_n} \varphi) \,\de\bm{x} \right) = 0.
 \end{equation*}
 Due to the fact that
 \begin{equation*}
   \lim_{n \rightarrow \infty} \int_{\Omega} f(\bm{x}) P_{\mathcal{D}_n} \varphi(\bm{x}) \,\de\bm{x}
   = \int_\Omega f(\bm{x}) \varphi(\bm{x}) \,\de\bm{x},
 \end{equation*}
 one gets that any limit $\underline{u}$ of a subsequence of solutions
 satisfies the weak problem \eqref{eq:isotropic_diffusion_weak_problem} with $v = \varphi$.
 Uniqueness of the solution to \eqref{eq:isotropic_diffusion_weak_problem} together with a classical density
 argument allow to deduce the convergence of the whole sequence $u$ to the weak problem solution
 $\underline{u}$ in $L^2(\Omega)$ as $h_{\mathcal{D}} \rightarrow 0$, since $\theta \leq \tilde\theta_{\mathcal{D}}$.
\end{proof*}

\section{Numerical results}
 \label{sec:numerical_results}
 
The \emph{Gauss corrected} finite volume diffusion scheme in Eq.\eqref{eq:face_gradient_corrected} was tested
on a number of different mesh types commonly adopted in industrial applications. These include meshes composed
of regular orthogonal hexahedra, skewed hexahedra, triangular prismatic and polygonal prismatic cells. All
these meshes were constructed by means of a commercial finite volume mesh generator \cite{Gambit}, which
usually produces meshes of acceptable non-orthogonality, as commonly required in practical applications. In
particular, the polygonal prismatic mesh was obtained after geometric dualization of the triangular prismatic
one. Notice also that finer meshes are not produced by conformal refinement techniques, but generated \emph{ex novo}.  It should be remarked again that the \emph{Gauss corrected} scheme reduces to the
two-point flux approximation on orthogonal meshes.

\begin{figure}[t!]
 \centering
 \includegraphics[width=0.35\linewidth]{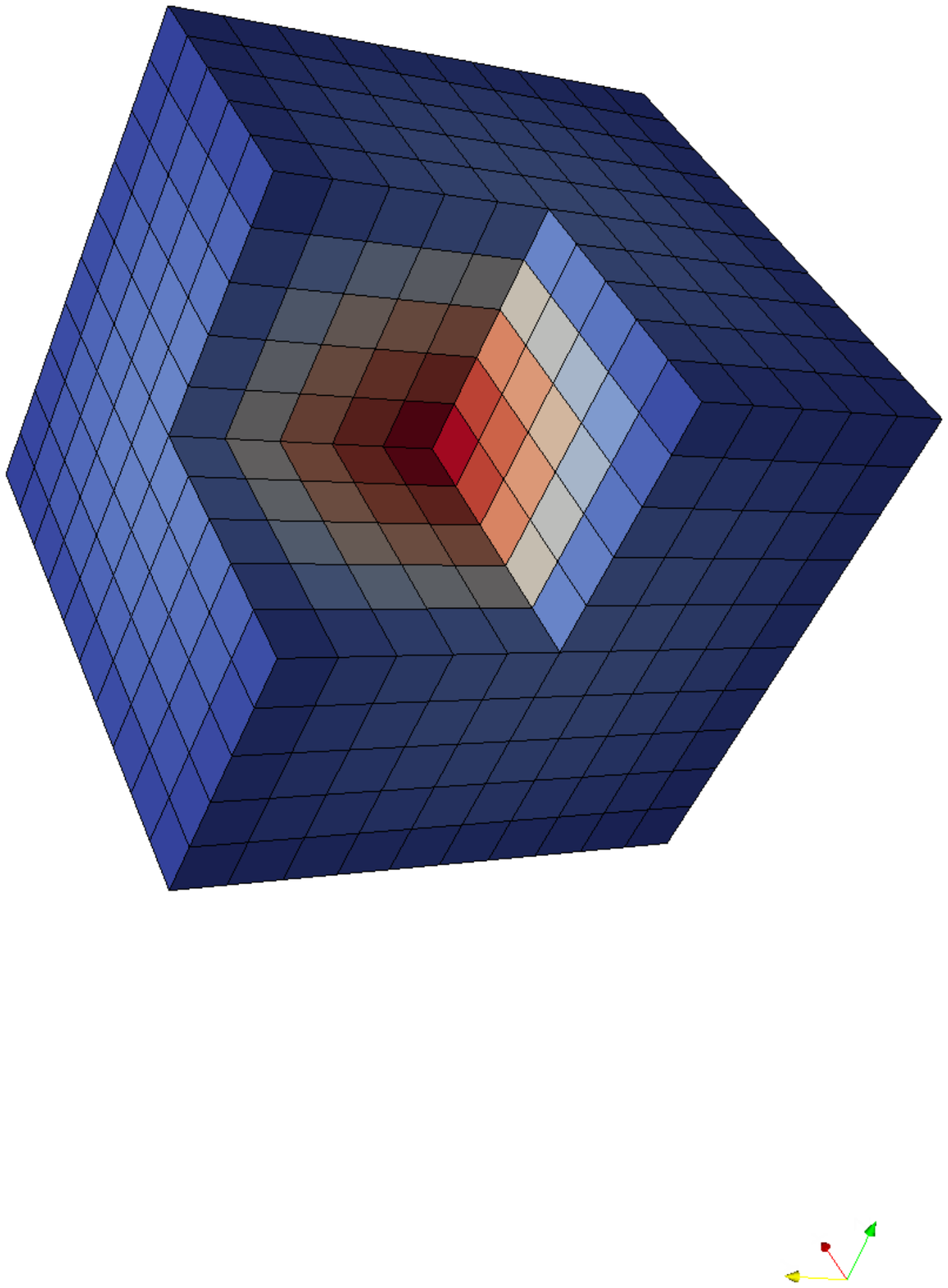} {$(a)$};
 \includegraphics[width=0.35\linewidth]{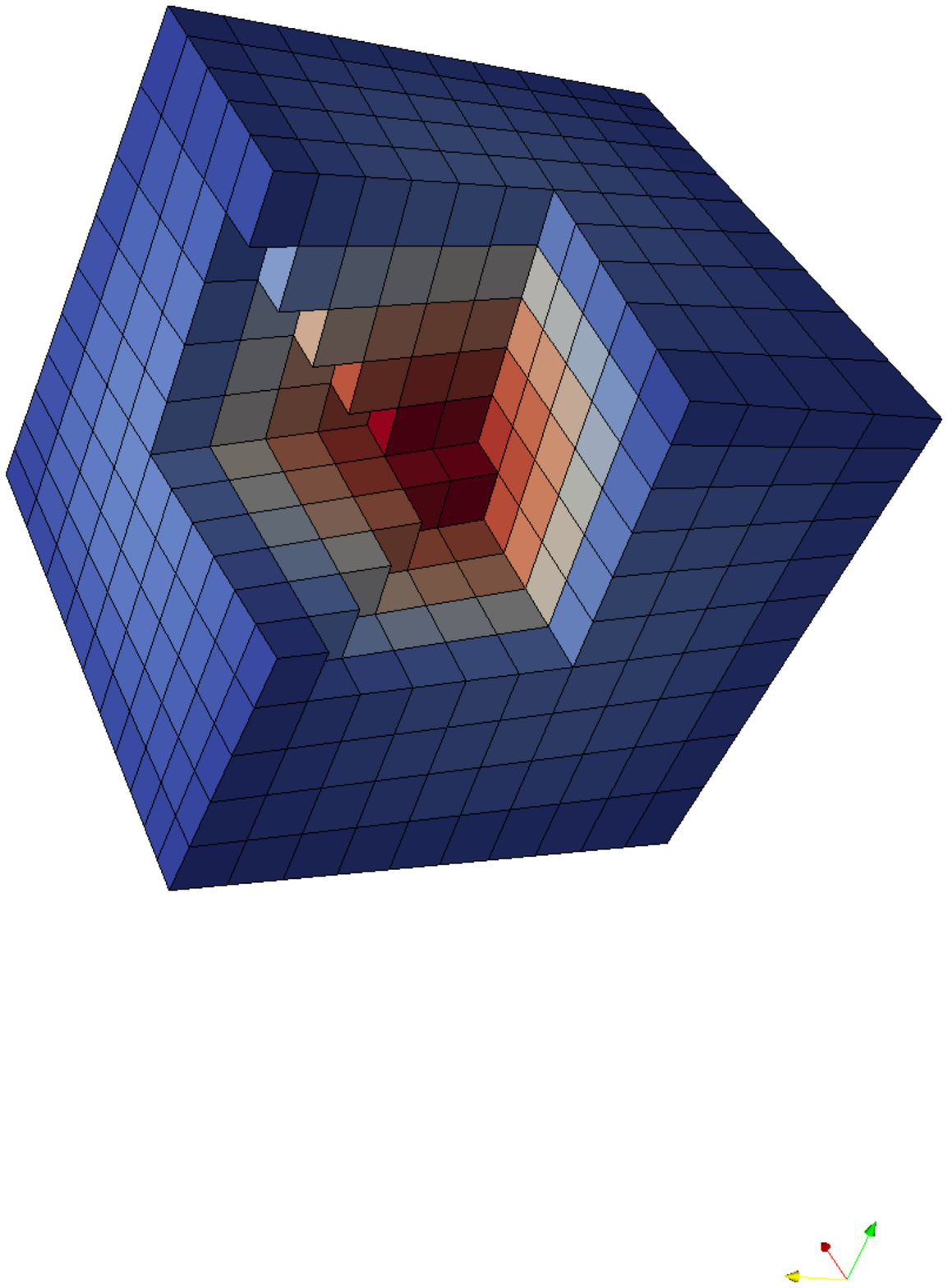}{$(b)$}
 \\
   \includegraphics[width=0.35\linewidth]{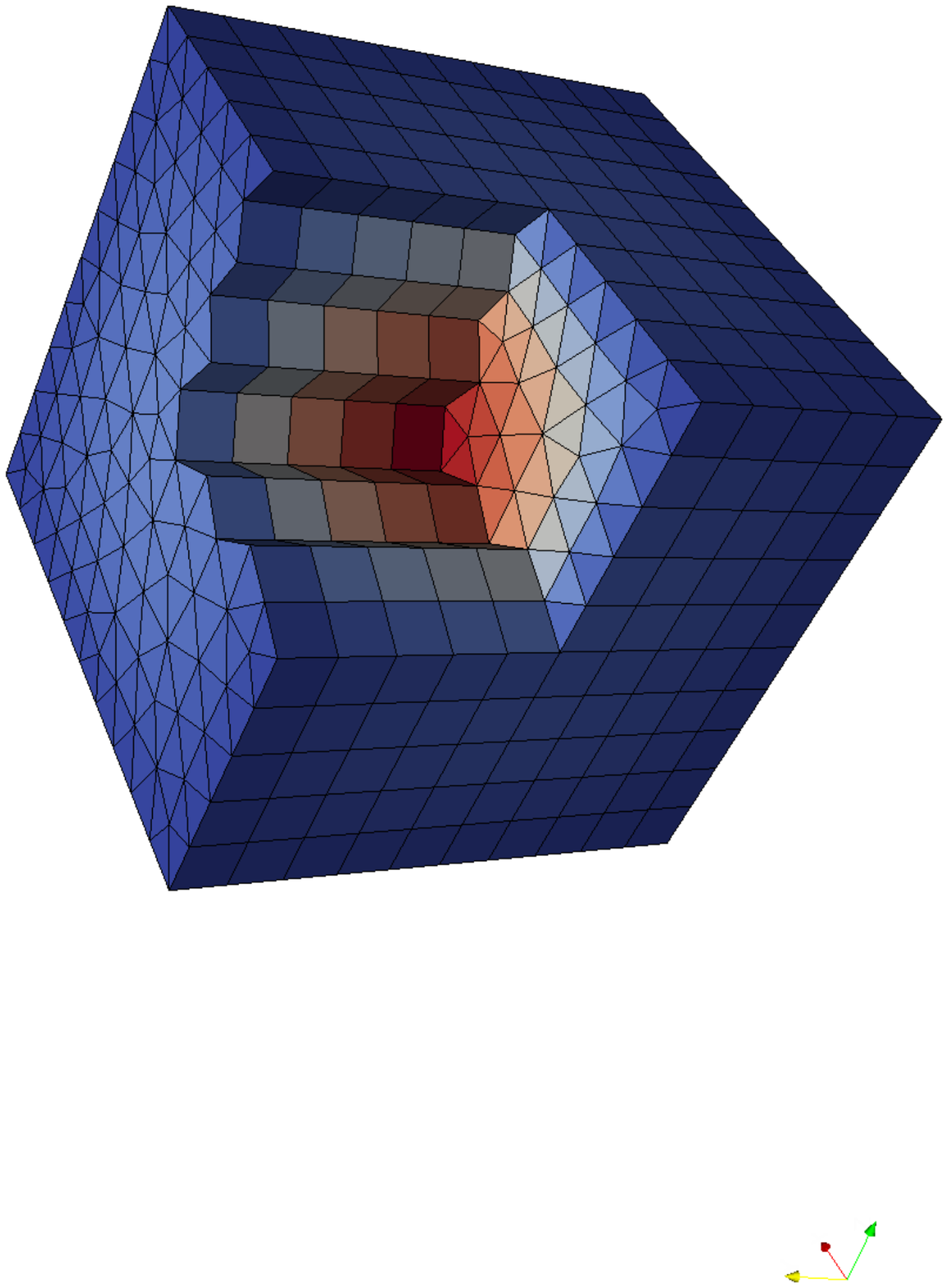} {$(c)$}
    \includegraphics[width=0.35\linewidth]{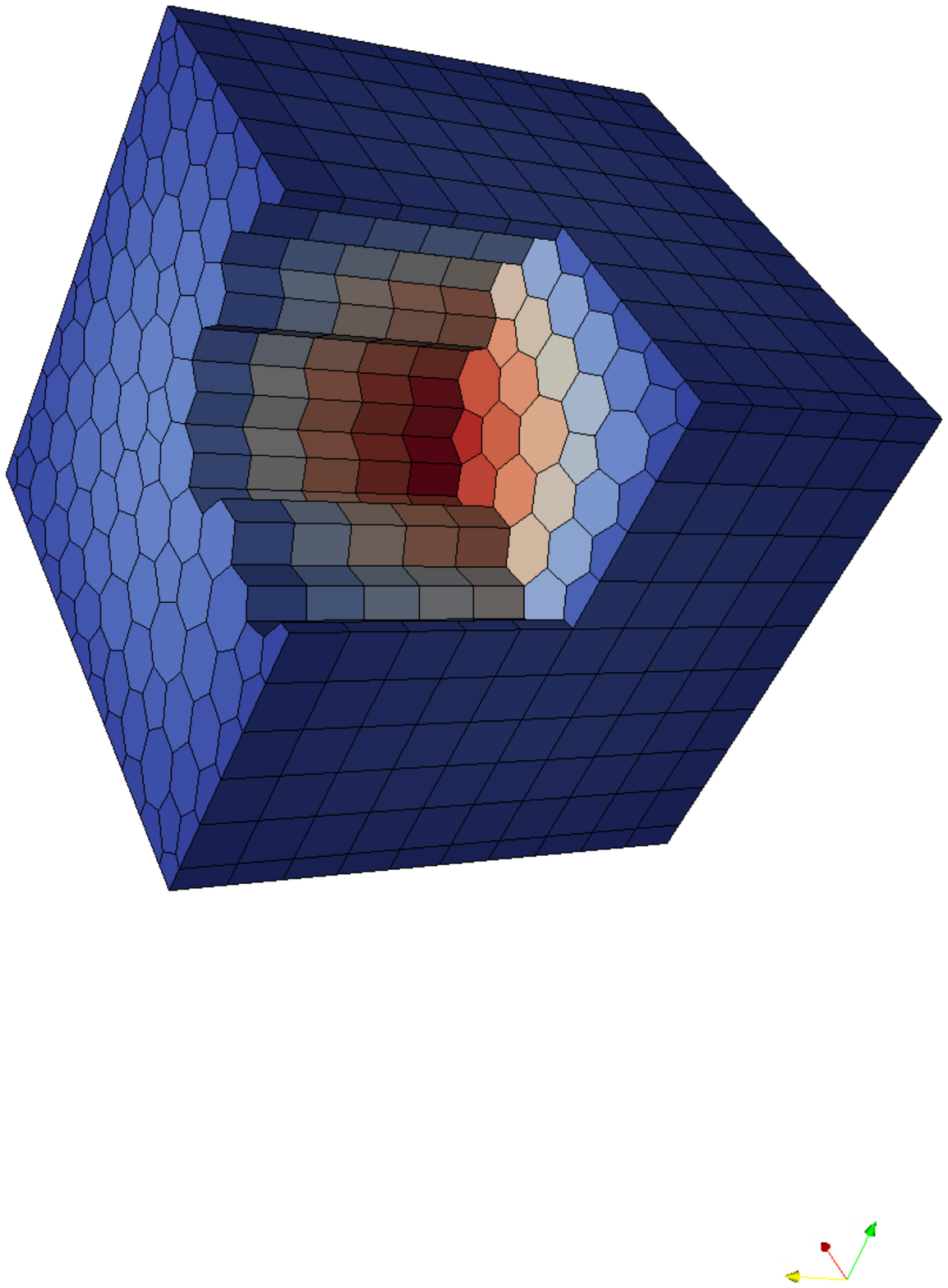} {$(d)$}
 \caption{Different mesh types used in the numerical test: ($a$) orthogonal hexahedral, ($b$) skewed
   hexahedral, ($c$) triangular prismatic, ($d$) polygonal prismatic.}
 \label{fig:meshes_test}
\end{figure}

\begin{table}[t!]
 \footnotesize
 \centering
 \setlength{\tabcolsep}{0.5em}
 \begin{tabularx}{0.71\textwidth}{l c c c c c}
 \toprule
  \emph{Mesh} & $\langle d \rangle_\Omega$ & $\langle \theta_\sigma \rangle$ & $\theta_\text{max}$
              & $\text{AR}_\text{max}$ & $\text{S}_{\text{max}}$ \\
 \midrule
hex       & $1.0000 \times 10^{-1}$ & 0      & 0      & 1     & 0 \\
          & $5.0000 \times 10^{-2}$ & 0      & 0      & 1     & 0 \\
          & $2.5000 \times 10^{-2}$ & 0      & 0      & 1     & 0 \\
          & $1.2500 \times 10^{-2}$ & 0      & 0      & 1     & 0 \\
hexSkew   & $1.0033 \times 10^{-1}$ & 9.017  & 15.138 & 2.085 & 0.145 \\
          & $5.0265 \times 10^{-2}$ & 8.982  & 15.392 & 2.352 & 0.162 \\
          & $2.5160 \times 10^{-2}$ & 8.953  & 15.531 & 2.520 & 0.175 \\
          & $1.2587 \times 10^{-2}$ & 8.935  & 15.588 & 2.615 & 0.181 \\
triPrism  & $7.4413 \times 10^{-2}$ & 3.396  & 11.059 & 2.963 & 0.221 \\
          & $3.7505 \times 10^{-2}$ & 2.862  & 10.651 & 2.969 & 0.218 \\
          & $1.8748 \times 10^{-2}$ & 2.780  & 10.652 & 3.146 & 0.241 \\
          & $9.3825 \times 10^{-3}$ & 2.574  & 11.008 & 3.176 & 0.239 \\
polyPrism & $9.7486 \times 10^{-2}$ & 5.246  & 14.912 & 2.790 & 0.741 \\
          & $4.9855 \times 10^{-2}$ & 3.945  & 14.587 & 2.778 & 0.739 \\
          & $2.5063 \times 10^{-2}$ & 3.416  & 15.082 & 2.777 & 0.739 \\
          & $1.2582 \times 10^{-2}$ & 3.022  & 14.802 & 2.777 & 0.739 \\
 \bottomrule
 \end{tabularx}
 \caption{Main geometric parameters for the different mesh types used in the accuracy test.}
  \label{tab:meshes}
\end{table}

The different mesh types are summarized in Table \ref{tab:meshes}, where relevant quantities are reported. 
These include the parameters normally observed as quality indices after the mesh generation process,
which are:
\begin{itemize}

\item \emph{Non-orthogonality}, measured by the angle between the line segment $(\bm{x}_L - \bm{x}_K)$,
 joining cell centroids adjacent to face $\sigma \in \mathcal{F}_K$, and the face normal $\bm{n}_{K,\sigma}$,
 that is
 \begin{eqnarray}
  \theta_\sigma &=& \arccos{\left(\frac{(\bm{x}_L - \bm{x}_K) \cdot \bm{n}_{K,\sigma}}{|\bm{x}_L - \bm{x}_K|}\right)},
    \nonumber \\
  \forall \sigma &\in& \mathcal{F}_{K,int}, \; \forall K \in \mathcal{M}.
 \end{eqnarray}
 A value close to $0$ is optimal , since it reduces the amount of non-orthogonal correction with respect to
 the two-point flux approximation, see, e.g., Eq.\eqref{eq:face_gradient_corrected}. Here, both the mean
 non-orthogonality angle $\langle \theta_\sigma \rangle = \frac{1}{\#\mathcal{F}_{int}} \sum_{\sigma \in
   \mathcal{F}_{int}} \theta_\sigma$ and the maximum non-orthogonality angle $\theta_\text{max} = \max_{\sigma
   \in \mathcal{F}_{int}} \theta_\sigma$ are considered.

\item \emph{Aspect ratio}, defined for each cell $K \in \mathcal{M}$ as
 \begin{equation}
  \text{AR}_K = \text{max} \left\{ \text{AR}_{BB(K)} , \frac{\sum_{i=1}^{d} |\sigma_{i}| }{6 |K|^{2/3}} \right\},
                 \quad \forall K \in \mathcal{M},
 \end{equation}
 where $\text{AR}_K$ is the bounding box aspect ratio
 \begin{equation}
  \text{AR}_{BB(K)} = \frac{\max_{\sigma_i \in \mathcal{F}_{BB(K)}} \left\{ |\sigma_i| \right\} }{\min_{\sigma_i \in \mathcal{F}_{BB(K)}} \left\{ |\sigma_i| \right\} }, \quad \forall K \in \mathcal{M},
 \end{equation}
 defined in terms of the cell bounding box $BB(K)$ which encloses the cell $K$ with a set of faces $\sigma_i$
 ($i=1,\ldots,d$) having normals oriented along the axes of the Cartesian reference frame used for the mesh
 definition. A value close to $1$ indicates that the cell is isotropic. Mesh statistics generally consider the
 maximum value of the cell aspect ratio $\text{AR}_\text{max} = \max_{K \in \mathcal{M}} \text{AR}_{K}$.

\item \emph{Skewness}, defined as the distance between the intersection point
 $\bm{y}_\sigma = [\bm{x}_L, \bm{x}_K] \cap \ol\sigma$ between the line segment $(\bm{x}_L - \bm{x}_K)$
 connecting adjacent cell centroids and separating face $\sigma \in \mathcal{F}_{int}$ and the face centroid
 $\bm{x}_\sigma$, that is
 \begin{equation}
  \text{S}_\sigma = \frac{|(\bm{x}_\sigma - \bm{y}_\sigma)|}{f_\sigma},
                    \quad \forall \sigma \in \mathcal{F},
 \end{equation}
 where the normalization factor $f_\sigma$ is
 \begin{subequations}
  \begin{align}
   &\begin{aligned}
     f_\sigma = \text{max} &\left\{ 0.2\, |(\bm{x}_L - \bm{x}_K)|, \right.\\
        &\hspace{1.0em}\left. \max_{i \in \mathcal{V}_\sigma} \left| (\bm{x}_i - \bm{x}_\sigma)
                    \cdot \frac{(\bm{x}_\sigma - \bm{y}_\sigma)}{|(\bm{x}_\sigma - \bm{y}_\sigma)|}  \right| \right\},
      \quad \forall \sigma \in \mathcal{F}_{int}
    \end{aligned}\\
   &\begin{aligned}
     f_\sigma = \text{max} &\left\{ 0.4\, |(\bm{y}_\sigma - \bm{x}_K)|, \right.\\
        &\hspace{1.0em}\left. \max_{i \in \mathcal{V}_\sigma} \left| (\bm{x}_i - \bm{x}_\sigma)
                    \cdot \frac{(\bm{x}_\sigma - \bm{y}_\sigma)}{|(\bm{x}_\sigma - \bm{y}_\sigma)|}  \right| \right\},
      \quad \forall \sigma \in \mathcal{F}_{ext}.
    \end{aligned}
  \end{align}
 \end{subequations}
 The optimal value for $S_\sigma$ is $0$, indicating that $\bm{y}_\sigma = \bm{x}_\sigma$, for which linear
 interpolation between adjacent cell values achieves second order consistency in face integral quantities.
 Mesh statistics generally take into account the maximum value of skewness $S_{\text{max}} = \max_{\sigma \in \mathcal{F}} S_\sigma$.
\end{itemize}
The mesh resolution is measured by the mean magnitude of the cell to cell distance, that is $\langle d
\rangle_\Omega = \langle | (\bm{x}_K - \bm{x}_L) | \rangle_{K|L \in \mathcal{F}_{int}}$.

\begin{figure}[t!]
\includegraphics[width=0.45\linewidth]{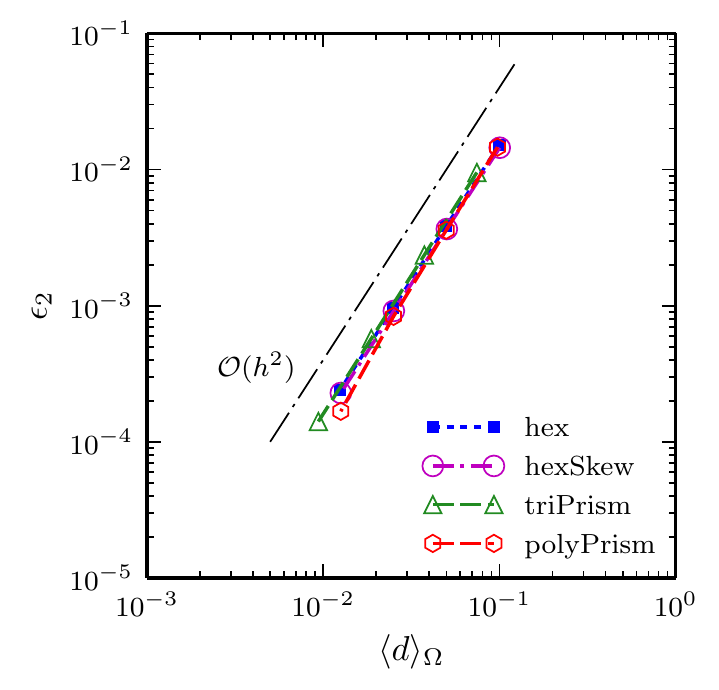}(a)
\includegraphics[width=0.45\linewidth]{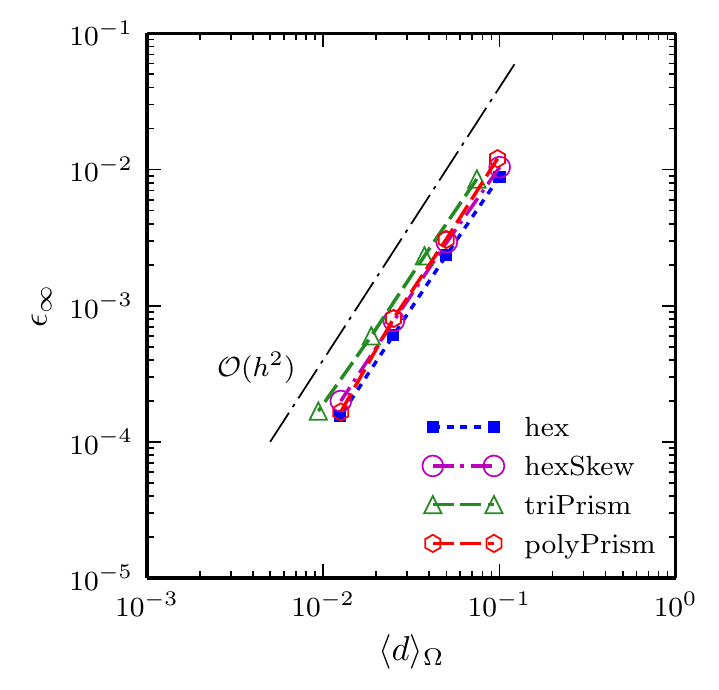}(b)
    \caption{Relative error curves for different mesh types: ($a$) $L^2$ relative error $\epsilon_2 = \| e \|_2 / \| \underline{u}
   \|_2 $, ($b$) $L^{\infty} $ relative error
   $\epsilon_2 = \| e \|_\infty / \| \underline{u} \|_\infty$.}
 \label{fig:errors}
\end{figure}

\begin{table}[t!]
 \footnotesize
 \centering
 \setlength{\tabcolsep}{0.5em}
 \begin{tabularx}{0.88\textwidth}{l c c c c c}
 \toprule
  \emph{Mesh} & $\langle d \rangle_\Omega$ & $\| e \|_{2}$ & $\| e \|_{\infty}$ & $p_2$ & $p_\infty$ \\
 \midrule
hex       & $1.0000 \times 10^{-1}$ & $9.2721 \times 10^{-5}$ & $1.3410 \times 10^{-4}$ & $1.984$ & $1.882$ \\
          & $5.0000 \times 10^{-2}$ & $2.3439 \times 10^{-5}$ & $3.6372 \times 10^{-5}$ & $1.996$ & $1.946$ \\
          & $2.5000 \times 10^{-2}$ & $5.8771 \times 10^{-6}$ & $9.4379 \times 10^{-6}$ & $1.999$ & $1.975$ \\
          & $1.2500 \times 10^{-2}$ & $1.4704 \times 10^{-6}$ & $2.4011 \times 10^{-6}$ & $-$ & $-$ \\
hexSkew   & $1.0033 \times 10^{-1}$ & $8.7858 \times 10^{-5}$ & $1.5863 \times 10^{-4}$ & $1.988$ & $1.809$ \\
          & $5.0265 \times 10^{-2}$ & $2.2232 \times 10^{-5}$ & $4.5431 \times 10^{-5}$ & $2.001$ & $1.914$ \\
          & $2.5160 \times 10^{-2}$ & $5.5665 \times 10^{-6}$ & $1.2077 \times 10^{-5}$ & $2.004$ & $1.957$ \\
          & $1.2587 \times 10^{-2}$ & $1.3890 \times 10^{-6}$ & $3.1135 \times 10^{-6}$ & $-$ & $-$ \\
triPrism  & $7.4413 \times 10^{-2}$ & $5.7868 \times 10^{-5}$ & $1.3136 \times 10^{-4}$ & $2.036$ & $1.886$ \\
          & $3.7505 \times 10^{-2}$ & $1.4340 \times 10^{-5}$ & $3.6076 \times 10^{-5}$ & $2.040$ & $1.946$ \\
          & $1.8748 \times 10^{-2}$ & $3.4852 \times 10^{-6}$ & $9.3566 \times 10^{-6}$ & $2.025$ & $1.835$ \\
          & $9.3825 \times 10^{-3}$ & $8.5803 \times 10^{-7}$ & $2.6272 \times 10^{-6}$ & $-$ & $-$ \\
polyPrism & $9.7486 \times 10^{-2}$ & $8.8853 \times 10^{-5}$ & $1.8380 \times 10^{-4}$ & $2.096$ & $2.011$ \\
          & $4.9855 \times 10^{-2}$ & $2.1796 \times 10^{-5}$ & $4.7717 \times 10^{-5}$ & $2.118$ & $1.941$ \\
          & $2.5063 \times 10^{-2}$ & $5.0781 \times 10^{-6}$ & $1.2557 \times 10^{-5}$ & $2.330$ & $2.292$ \\
          & $1.2582 \times 10^{-2}$ & $1.0198 \times 10^{-6}$ & $2.5873 \times 10^{-6}$ & $-$ & $-$ \\
 \bottomrule
 \end{tabularx}
 \caption{Error behaviour in the numerical convergence test.}
  \label{tab:errors}
\end{table}

%
%

Numerical experiments were carried out on $\Omega = (0,1) \times (0,1) \times (0,1)$ assuming $\alpha(\bm{x}) = 1$
and considering the exact solution of problem \eqref{eq:isotropic_diffusion_problem} given by
$\underline{u}(x_1, x_2, x_3) = x_1(1 - x_1) x_2(1 - x_2) x_3(1 - x_3)$. On each mesh, the error is measured
as $e(\bm{x}_K) = u_K - \underline{u}(\bm{x}_K)$ for $K \in \mathcal{M}$ and it allows to estimate
empirically the rate of convergence between two successive mesh sizes. The \emph{Gauss corrected} scheme is
implemented with the deferred correction approach, with the non-orthogonal correction term implemented
explicitly and thus requiring outer iterations which are terminated with a tolerance level of $1.0 \times
10^{-4}$. The associated linear system is solved by a preconditioned conjugate gradient method with
tolerance $1.0 \times 10^{-16}$ and DIC preconditioning.

The error norms $\| e \|_2$ and $\| e \|_\infty$ together with the corresponding empirical orders of
convergence $p_2$ and $p_\infty$ are reported in Table \ref{tab:errors}, while the relative errors norms
$\epsilon_2 = \| e \|_2 / \| \underline{u} \|_2 $ and $\epsilon_2 = \| e \|_\infty / \| \underline{u}
\|_\infty$ are shown in Figure \ref{fig:errors}. From both quantities, it is evident that second order
accuracy is empirically verified for hexahedral, skewed hexahedral, triangular prismatic and polyhedral
prismatic mesh types. It is remarkable that the accuracy of the \emph{Gauss corrected} scheme appears insensitive to the cell shape, with only minimal differences in the infinity norm.

\section{Beyond the  Gauss gradient scheme}
 \label{sec:least_squares_tests} \indent

The  Gauss discrete gradient operator $\nabla_{\mathcal{D}}$
 that was introduced in Eq.\eqref{eq:discrete_gradient} is
bounded, weakly convergent and consistent, but generally it is not coercive. Thus, on strongly
non-orthogonal meshes, the correction term in the \emph{Gauss corrected} approach may not be coercive and
consequently   hamper the convergence of the finite volume scheme.

\begin{figure}[t!]
 \centering
\includegraphics[width=0.35\linewidth]{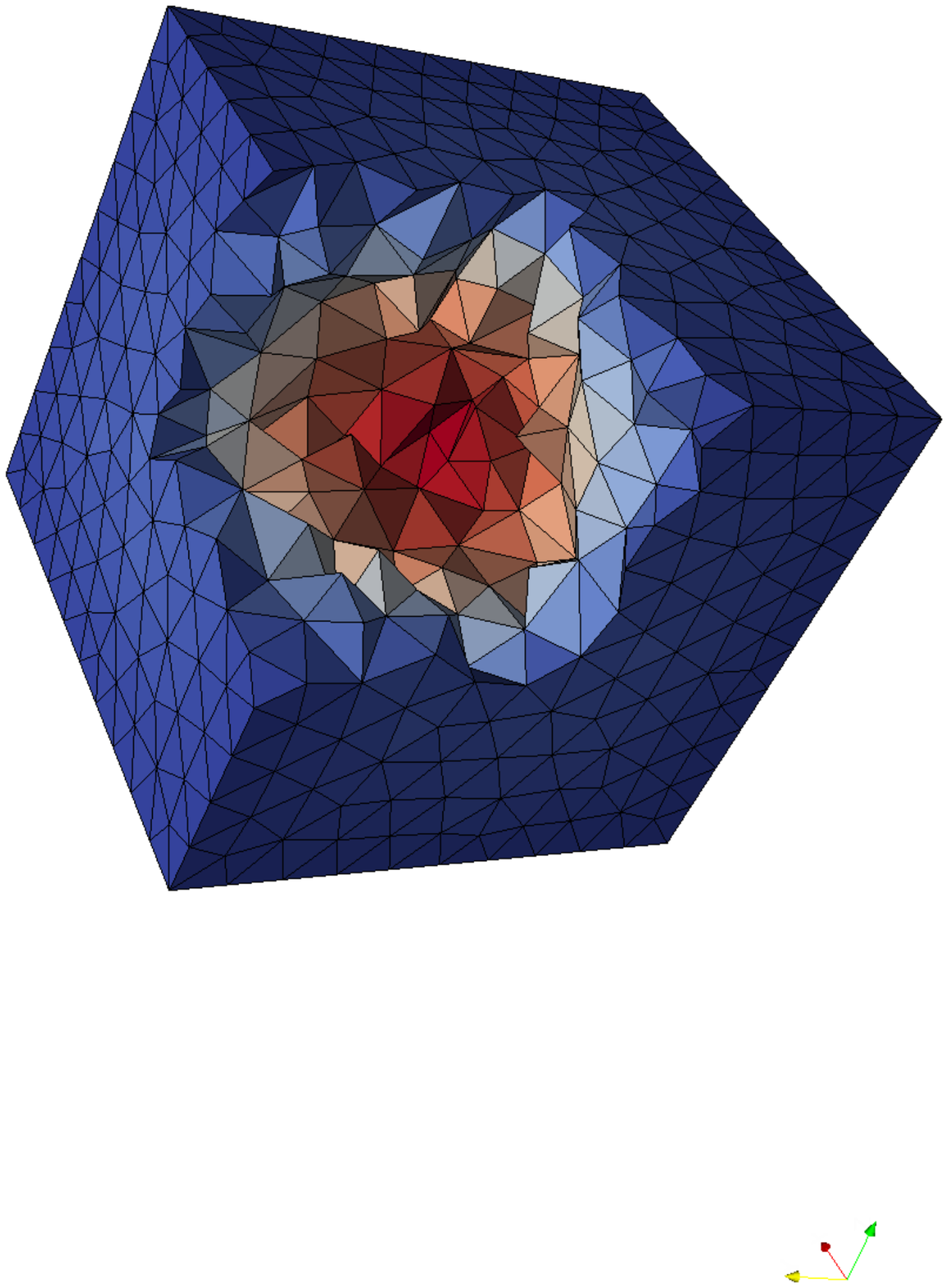} {$(a)$}
 \includegraphics[width=0.35\linewidth]{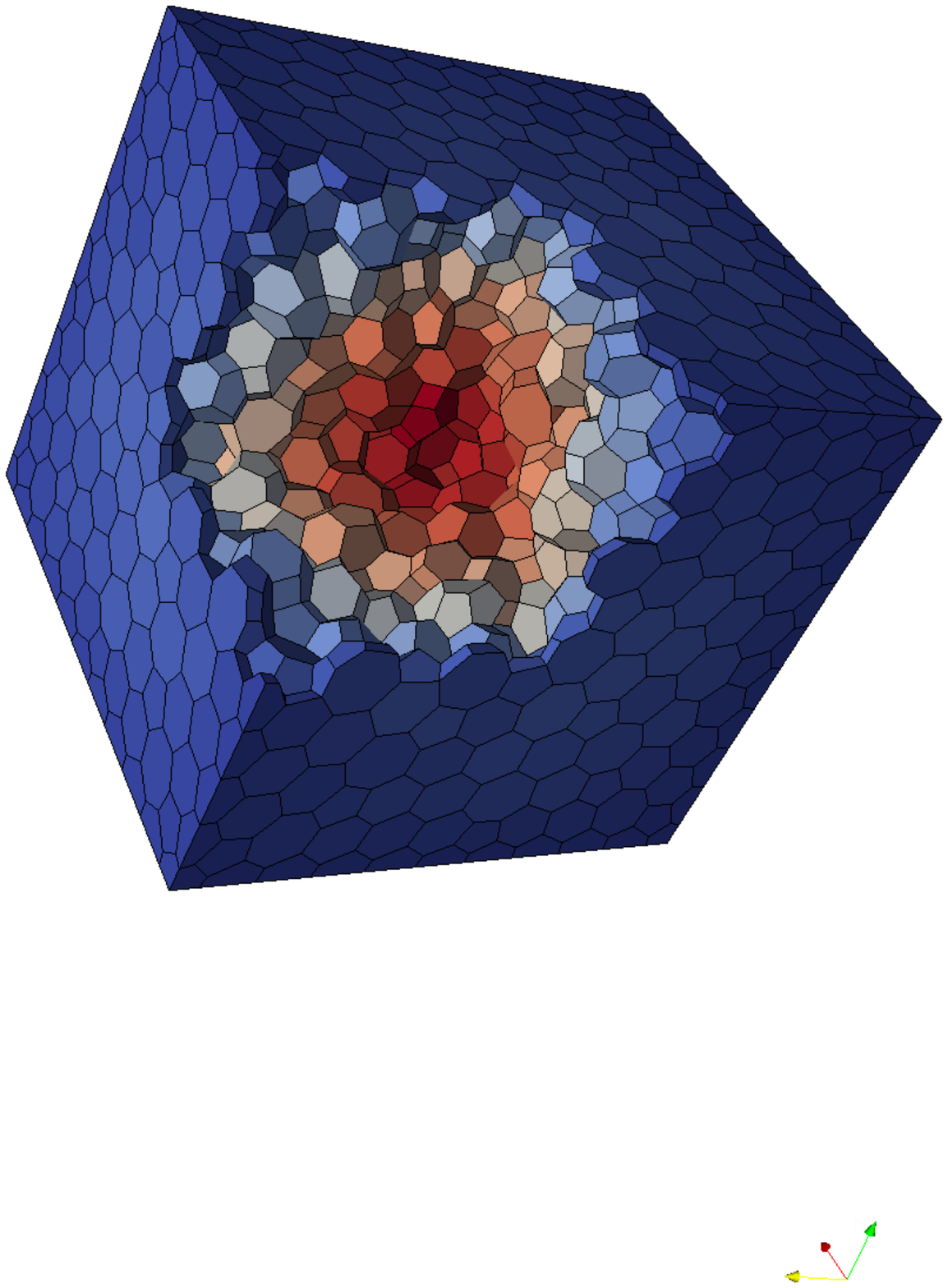} {$(b)$}
 \caption{Strongly non-orthogonal mesh types used in the numerical test: ($a$) tetrahedral and ($b$)
   polyhedral meshes.}
 \label{fig:meshes_test_non_orthogonal}
\end{figure}

To verify this point empirically, the same diffusion problem studied empirically in section
\ref{sec:numerical_results} is now solved on a sequence of highly
non-orthogonal tetrahedral and polyhedral meshes, see Figure \ref{fig:meshes_test_non_orthogonal}, whose
geometric parameters are summarized in Table \ref{tab:meshes_non_orthogonal}. It is important to notice that
the maximum non-orthogonality angle is such thats $\theta_\text{max} > \pi/4$ almost for every mesh, with the only exception of
the two coarsest polyhedral meshes. This implies that the non-orthogonal correction term is the dominant term
in the numerical flux.

The \emph{Gauss corrected} scheme can still be applied, provided that a coercive gradient scheme is adopted
for the non-orthogonal correction term. To this end, we assess here the performance a gradient approximation based on a  least square fit, based on the fact that in   linear upwind schemes it is known empirically   to provide a coercive
gradient discretization in the case of highly non-orthogonal tetrahedral meshes. Notice that, on orthogonal
meshes it reduces to the  Gauss scheme, hence becoming non-coercive. But this is of no concern as long
as it is adopted only for the construction of the non-orthogonal correction term in the \emph{Gauss corrected}
fluxes.

\begin{table}[t!]
 \footnotesize
 \centering
 \setlength{\tabcolsep}{0.5em}
 \begin{tabularx}{0.71\textwidth}{l c c c c c}
 \toprule
  \emph{Mesh} & $\langle d \rangle_\Omega$ & $\langle \theta_\sigma \rangle$ & $\theta_\text{max}$
              & $\text{AR}_\text{max}$ & $\text{S}_{\text{max}}$ \\
 \midrule
tet       & $4.2582 \times 10^{-2}$ & 19.058 & 59.778 & 5.614 & 0.569 \\
          & $2.1113 \times 10^{-2}$ & 19.307 & 63.201 & 7.166 & 0.743 \\
          & $1.1023 \times 10^{-2}$ & 19.597 & 65.714 & 8.468 & 0.926 \\
          & $5.5488 \times 10^{-3}$ & 19.741 & 66.510 & 8.370 & 0.901 \\
poly      & $1.0481 \times 10^{-1}$ & 11.688 & 38.598 & 4.592 & 1.062 \\
          & $5.3575 \times 10^{-2}$ & 11.801 & 40.406 & 3.188 & 1.134 \\
          & $2.8269 \times 10^{-2}$ & 11.925 & 47.784 & 4.048 & 1.234 \\
          & $1.4341 \times 10^{-2}$ & 11.944 & 50.781 & 4.188 & 1.439 \\
 \bottomrule
 \end{tabularx}
 \caption{Main geometric parameters for the highly non-orthogonal mesh types used in the accuracy test.}
  \label{tab:meshes_non_orthogonal}
\end{table}

Following \cite{barth:2004}, on a non-orthogonal mesh like that of Definition \ref{def:polyhedral_mesh}, we define
the discrete gradient operator $\nabla_{\mathcal{D}}^{LS} : H_\mathcal{D}(\Omega) \rightarrow H_\mathcal{D}(\Omega)^d$
 as the piecewise constant function
\begin{equation}
  \nabla_K^{LS} u = \sum_{\sigma \in \mathcal{F}_{K,int}} \left( u_L - u_K \right) \bm{v}_{K,\sigma}
                  + \sum_{\sigma \in \mathcal{F}_{K,ext}} \left( u_\sigma - u_K \right) \bm{v}_{K,\sigma}
 \label{eq:least_squares_gradient}
\end{equation}
for $ u \in H_\mathcal{D}(\Omega)$, where the least squares vectors
\begin{subequations}
 \begin{align}
   \bm{v}_{K,\sigma} &= w_{K,\sigma}\, \bm{W}_K^{-1} \left( \bm{x}_L - \bm{x}_K \right),
                \;\forall \sigma \!\in\! \mathcal{F}_{K,int}\\
   \bm{v}_{K,\sigma} &= w_{K,\sigma}\, \bm{W}_K^{-1} \left( \bm{x}_\sigma - \bm{x}_K \right),
                \;\forall \sigma \!\in\! \mathcal{F}_{K,ext}
 \end{align}
 \label{eq:least_squares_vectors}%
\end{subequations}
are defined from the weighting tensor
\begin{equation}
 \begin{aligned}
  \bm{W}_K =
  &\sum_{\sigma \in \mathcal{F}_{K,int}} w_{K,\sigma} \left( \bm{x}_L - \bm{x}_K \right) \left( \bm{x}_L - \bm{x}_K \right)^\intercal\\
  &+ \sum_{\sigma \in \mathcal{F}_{K,ext}} w_{K,\sigma} \left( \bm{x}_\sigma - \bm{x}_K \right) \left( \bm{x}_\sigma - \bm{x}_K \right)^\intercal,
  \;\forall K \!\in\! \mathcal{M}
 \end{aligned}
 \label{eq:weighting_tensor}
\end{equation}
as well as from the face weights $w_{K,\sigma}$, $\forall \sigma \in \mathcal{F}_K$ and $\forall K \in \mathcal{M}$.
Different expressions for the face weights can be adopted. Here, we use the formulae
\begin{subequations}
 \begin{align}
   w_{K,\sigma} &= \frac{d_{K,\sigma}}{d_{K,L}} \frac{|\sigma|}{|\left( \bm{x}_L - \bm{x}_K \right)|^2},
                \;\forall \sigma \!\in\! \mathcal{F}_{K,int}\\
   w_{K,\sigma} &= \frac{|\sigma|}{|\left( \bm{x}_L - \bm{x}_K \right)|^2},
                \;\forall \sigma \!\in\! \mathcal{F}_{K,ext} .
 \end{align}
 \label{eq:least_squares_weights}%
\end{subequations}
The  least squaress gradient scheme is empirically constructed from the approximate Taylor expansion at
adjacent cell centroids and face centroids
\begin{subequations}
 \begin{align}
  u (\bm{x}_L) &\approx u_K + \nabla_K^{LS} u \cdot \left( \bm{x}_L - \bm{x}_K \right),
                       \;\forall \sigma \!\in\! \mathcal{F}_{K,int},\\
  u (\bm{x}_\sigma) &\approx u_K + \nabla_K^{LS} u \cdot \left( \bm{x}_\sigma - \bm{x}_K \right),
                       \;\forall \sigma \!\in\! \mathcal{F}_{K,ext},
 \end{align}
 \label{eq:least_squares_construction}%
\end{subequations}
requiring the minimization of the piecewise constant mean-square-error objective function
\begin{equation*}
 \begin{aligned}
  G_K = &\sum_{\sigma \in \mathcal{F}_{K,int}} w_{K,\sigma} \left( u_L - u_K - \nabla_K^{LS} u \cdot \left( \bm{x}_L - \bm{x}_K \right) \right)^2\\
       &+ \sum_{\sigma \in \mathcal{F}_{K,ext}} w_{K,\sigma} \left( u_\sigma - u_K - \nabla_K^{LS} u \cdot \left( \bm{x}_\sigma - \bm{x}_K \right) \right)^2.
 \end{aligned}
\end{equation*}
%

\begin{figure}[t!]
\includegraphics[width=50mm]{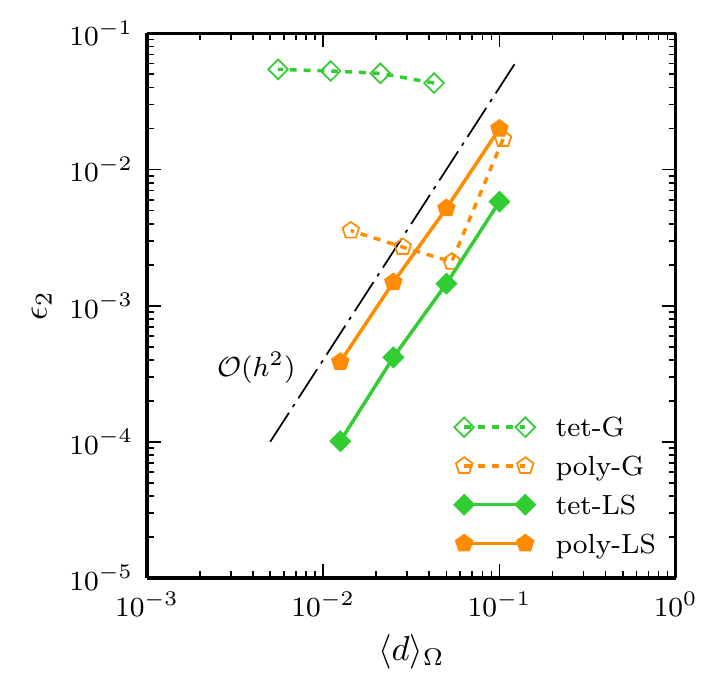}(a)
\includegraphics[width=50mm]{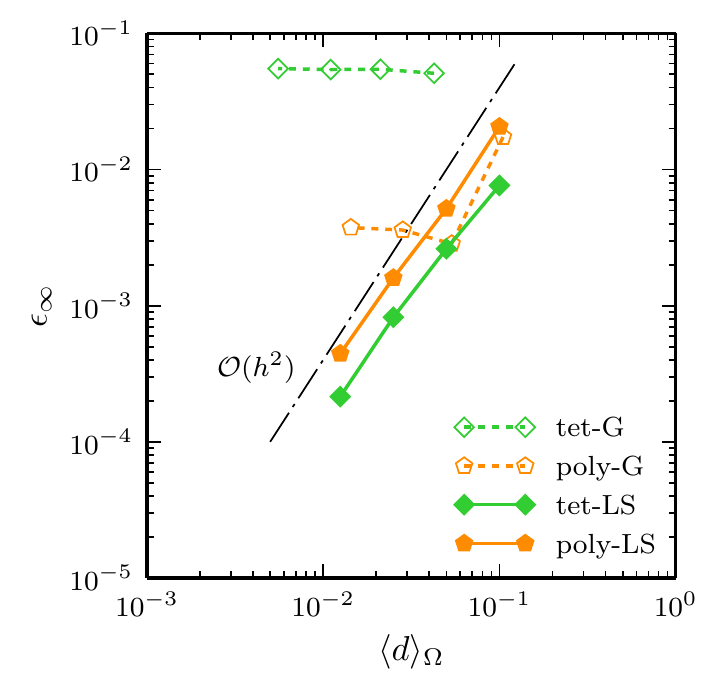}(b)
    \caption{Relative error curves for non-orthgonal mesh types: ($a$) $L^2$ relative error $\epsilon_2 = \| e
   \|_2 / \| \underline{u} \|_2 $, ($b$) $L^{\infty} $ relative error $\epsilon_2 = \| e \|_\infty / \|
   \underline{u} \|_\infty$. The \emph{Gauss} (G) gradient scheme is assessed against the \emph{leastSquares}
   (LS) scheme.}
 \label{fig:errors_non_orthogonal}
\end{figure}

\begin{table}[t!]
 \footnotesize
 \centering
 \setlength{\tabcolsep}{0.5em}
 \begin{tabularx}{0.94\textwidth}{l c c c c c c}
 \toprule
  \emph{Mesh} & \emph{Grad} & $\langle d \rangle_\Omega$ & $\| e \|_{2}$ & $\| e \|_{\infty}$ & $p_2$ & $p_\infty$ \\
 \midrule
tet       & G & $4.2582 \times 10^{-2}$ & $2.6231 \times 10^{-4}$ & $7.9063 \times 10^{-4}$ & $-0.234$ & $-0.104$ \\
          & G & $2.1113 \times 10^{-2}$ & $3.0907 \times 10^{-4}$ & $8.5032 \times 10^{-4}$ & $-0.062$ & $0.007$ \\
          & G & $1.1023 \times 10^{-2}$ & $3.2172 \times 10^{-4}$ & $8.4668 \times 10^{-4}$ & $-0.039$ & $-0.022$ \\
          & G & $5.5488 \times 10^{-3}$ & $3.3056 \times 10^{-4}$ & $8.5973 \times 10^{-4}$ & $-$ & $-$ \\
tet       & LS & $4.2582 \times 10^{-2}$ & $3.5338 \times 10^{-5}$ & $1.1860 \times 10^{-4}$ & $1.978$ & $1.518$ \\
          & LS & $2.1113 \times 10^{-2}$ & $8.8216 \times 10^{-6}$ & $4.0872 \times 10^{-5}$ & $1.917$ & $1.782$ \\
          & LS & $1.1023 \times 10^{-2}$ & $2.5376 \times 10^{-6}$ & $1.2836 \times 10^{-5}$ & $2.064$ & $1.955$ \\
          & LS & $5.5488 \times 10^{-3}$ & $6.1543 \times 10^{-7}$ & $3.3552 \times 10^{-6}$ & $-$ & $-$ \\
poly      & G & $1.0481 \times 10^{-1}$ & $1.0196 \times 10^{-4}$ & $2.6907 \times 10^{-4}$ & $3.094$ & $2.685$ \\
          & G & $5.3575 \times 10^{-2}$ & $1.2782 \times 10^{-5}$ & $4.4394 \times 10^{-5}$ & $-0.388$ & $-0.366$ \\
          & G & $2.8269 \times 10^{-2}$ & $1.6381 \times 10^{-5}$ & $5.6089 \times 10^{-5}$ & $-0.411$ & $-0.061$ \\
          & G & $1.4341 \times 10^{-2}$ & $2.1646 \times 10^{-5}$ & $5.8460 \times 10^{-5}$ & $-$ & $-$ \\
poly      & LS & $1.0481 \times 10^{-1}$ & $1.2148 \times 10^{-4}$ & $3.2076 \times 10^{-4}$ & $2.009$ & $2.061$ \\
          & LS & $5.3575 \times 10^{-2}$ & $3.1556 \times 10^{-5}$ & $8.0425 \times 10^{-5}$ & $1.957$ & $1.830$ \\
          & LS & $2.8269 \times 10^{-2}$ & $9.0307 \times 10^{-6}$ & $2.4964 \times 10^{-5}$ & $1.983$ & $1.886$ \\
          & LS & $1.4341 \times 10^{-2}$ & $2.3513 \times 10^{-6}$ & $6.9417 \times 10^{-6}$ & $-$ & $-$ \\
 \bottomrule
 \end{tabularx}
 \caption{Error behaviour for \emph{Gauss} (G) and \emph{leastSquares} (LS) gradient schemes in the numerical
   convergence test on highly non-orthogonal meshes.}
  \label{tab:errors_non_orthogonal}
\end{table}

The results of the comparison between the  Gauss and  least squares gradient schemes in the
construction of the non-orthogonal term inside the \emph{Gauss corrected} diffusion scheme are reported in
Figure \ref{fig:errors_non_orthogonal} and Table \ref{tab:errors_non_orthogonal}. On the strongly
non-orthogonal tetrahedral meshes, the  Gauss gradient does not allow to obtain a coercive numerical
flux and leads to stagnation. The same discrete gradient scheme converges only on
the first two coarser polyhedral meshes, while it diverges again on the two finest polyhedral meshes.  On the
contrary, when the non-orthogonal correction term of the \emph{Guss corrected} scheme is constructed from the
 least squares gradient scheme,   a convergent
second order behaviour is recovered on all the highly non-orthogonal meshes considered here. 
Even though in the present
investigation no analytical results have been obtained for this discrete gradient scheme, it appears to be
able to overcome the main limitation of the Gauss discrete gradient on strongly non-orthogonal meshes.

\section{Conclusions}
 \label{sec:conclusions}

In this work, we have  proven the convergence of the \emph{Gauss corrected} scheme on
unstructured meshes satisfying a global and rather weak mesh regularity condition. 
 This goal has been achieved adapting  the approach used in \cite{eymard:2006b} for the
convergence analysis of a cell-centered finite volume scheme for anisotropic diffusion problems on orthogonal
meshes. We have also shown empirically that a least square approach to the gradient computation  can
  provide second order convergence even when the mild mesh regularity condition is violated.
To the best of the authors' knowledge, the convergence properties of the finite volume method 
analyzed here have never been
studied rigorously in the case of non-orthogonal meshes. 
 Indeed, convergence analyses of finite volume schemes for
diffusion operators on unstructured mesh types are usually limited to polyhedral meshes satisfying an
orthogonality condition \cite{eymard:2000}, \cite{eymard:2006b}. This is quite restrictive in practice, since
none of the robust mesh generators usually adopted for pre-processing of industrial configurations are able
to guarantee this condition.

From our development, it can be seen how the analysis of finite volume schemes is greatly simplified if it is approached from the associated discrete
weak formulation using the functional tools defined in \cite{eymard:2007d}. In particular, from this approach
several interesting conclusions can be drawn, without recourse to the classical consistency analysis in terms
of Taylor series expansion. 
In the case of industrial finite volume schemes, such as the one analyzed here, these conclusions are
particularly interesting, because they shed some light over the properties of techniques
for which  typically only empirical results are available.
The role of the discrete gradient scheme is fundamental in many terms of the associated weak form, and its
relevance is also reflected in the finite volume formulation, even though it may not be completely apparent
when starting directly from the flux balance equations.  
By working with the weak form,  conclusions about  the coercivity of the finite volume scheme can be drawn,
which are often considered inaccessible in the finite volume framework.
Relevant discrete functional analysis results can be applied directly, in order to establish conditions for
convergence and eventually to obtain error estimates under sufficient regularity assumptions.
In particular, this approach allows to identify the mesh regularity requirements and the sufficient conditions for
convergence, as well as to suggest possible future improvements.

 The basic idea of correcting the two-point flux
approximation with an additional term accounting for the local mesh non - orthogonality can also be found in
many other finite volume schemes \cite{coudiere:1999}, \cite{lai:2000}, \cite{mathur:1997},
\cite{moukalled:2016}, \cite{moureau:2011}, \cite{muzaferija:1997}, \cite{perezsegarra:2006a}, and thus
similar analysis techniques could also be applied to investigate sufficient conditions for convergence of other
diffusion schemes. Connections with the asymmetric gradient schemes recently proposed in
\cite{droniou:2017} also suggest a possible alternative to the present analysis.
Finally, the diffusion operator analyzed in this work  was used in \cite{dellarocca:2018} as the basis for the construction of accurate and efficient stabilized pressure
correction methods for colocated finite volume schemes. The properties and advantages of these methods will be discussed in a series of forthcoming companion papers.


\section*{Acknowledgements}
The results presented in this work are part of the doctoral thesis in Applied Mathematics
\cite{dellarocca:2018} discussed by A.D.R. at Politecnico di Milano in 2018. The comments on the thesis by
 D. Di Pietro and J. Szmelter  are kindly acknowledged. A.D.R. would also like to thank Tenova S.p.A. for
sponsoring his Executive PhD at Politecnico di Milano and all the faculty members at MOX for their
support.\\


\section*{Appendix: Properties of the Gauss gradient scheme}
 \label{sec:appendix} \indent

The Gauss gradient operator $\nabla_{\mathcal{D}}$ introduced in Eq.\eqref{eq:discrete_gradient} is
bounded, weakly convergent and consistent.  In order to prove these properties, it is sufficient to rewrite it in the form
\begin{equation}
  \nabla_K u = \frac{1}{|K|} \sum_{\sigma \in \mathcal{F}_K} |\sigma|
                                             \frac{d_{K,\sigma}}{d_{K,L}}\left( u_L - u_K \right) \bm{n}_{K,\sigma}
 \label{eq:discrete_gradient_gauss}
\end{equation}
which allows to prove that is is bounded in the $L^2(\Omega)^d$-norm.

\begin{lemma}[Bound on $\nabla_{\mathcal{D}} u$]
 Let $\Omega$ be a bounded open connected polyhedral subset of $\mathbb{R}^d$, $d \in \mathbb{N}^\star$. Let
 $\mathcal{D}$ be an admissible finite volume discretization in sense of Definition \ref{def:polyhedral_mesh}
 and let $0 < \theta \leq \tilde\theta_\mathcal{D}$. Then, the exists $C$ depending only on $d$, $\alpha$ and
 $\theta$ such that, for all $u \in H_\mathcal{D}$, one has
 \begin{equation}
   \| \nabla_{\mathcal{D}} u \|_{L^2(\Omega)^d} \leq C \, \|u\|_\mathcal{D}.
  \label{eq:gauss_gradient_bound}
 \end{equation}
 \label{lem:gauss_gradient_bound}
\end{lemma}

\begin{proof*}
 Let $u \in H_\mathcal{D}$. As in Lemma \ref{lem:directional_gradient_bound}, one introduces, for
 all $K \in \mathcal{M}$, $L \in \mathcal{N}_K$ and $\sigma = K|L,$ the difference quantities
 $\delta_{K,\sigma} \bm{x}$ and $\delta_{K,\sigma} u$, from which the inner product norm in
 \eqref{eq:parallel_inner_product} leads for a given $K \in \mathcal{M}$ to
 \begin{equation*}
   \|u\|_{\mathcal{D}}^2
        = \sum_{K \in \mathcal{M}} \frac{1}{2} \sum_{L \in \mathcal{N}_K} \tau_{K|L} (\delta_{K,L} u)^2
          + \sum_{K \in \mathcal{M}} \sum_{\sigma \in \mathcal{F}_{K,ext}} \tau_{K,\sigma} (\delta_{K,\sigma} u)^2.
 \end{equation*}
 Then, from Eq.\eqref{eq:discrete_gradient_gauss} one obtains that
 \begin{equation*}
   |K| (\nabla u)_K
   = \sum_{\sigma \in \mathcal{F}_{K}} \tau_\sigma\, d_{K,\sigma} \delta_{K,\sigma} u.
 \end{equation*}
 By using the Cauchy-Schwartz inequality, one obtains that
 \begin{equation*}
  \begin{aligned}
   |K|^2 \, \left| (\nabla u)_K \right|^2 &\leq
   \sum_{\sigma \in \mathcal{F}_K } \tau_\sigma \left| d_{K,\sigma} \bm{n}_{K,\sigma} \right|^2 
      \sum_{\sigma \in \mathcal{F}_K } \tau_\sigma \left( \delta_{K,\sigma} u \right)^2 \\
   &\leq \sum_{\sigma \in \mathcal{F}_K } d \left| D_{K,\sigma} \right| \tau_\sigma
     \sum_{\sigma \in \mathcal{F}_K } \tau_\sigma \left( \delta_{K,\sigma} u \right)^2 \\
   &= d \, |K|\, \sum_{\sigma \in \mathcal{F}_K } \tau_\sigma \left( \delta_{K,\sigma} u \right)^2 .
  \end{aligned}
 \end{equation*}
 Summing over all $K \in \mathcal{M}$, one obtains that
 \begin{equation*}
  \begin{aligned}
   &\sum_{K \in \mathcal{M}} |K| \left| \nabla_K u \right|^2
   \leq d \sum_{K \in \mathcal{M}} \sum_{\sigma \in \mathcal{F}_K } \tau_\sigma \left( \delta_{K,\sigma} u \right)^2\\
   &\hspace{4.0em}\leq d \sum_{K \in \mathcal{M}} 
     \left( \sum_{\sigma \in \mathcal{F}_{K,int}} \tau_\sigma \left( \delta_{K,\sigma} u \right)^2 
     + 2 \sum_{\sigma \in \mathcal{F}_{K.ext}} \tau_\sigma \left( \delta_{K,\sigma} u \right)^2 \right)\\
   &\hspace{4.0em}= 2\, d \|u\|_{\mathcal{D}}^2
  \end{aligned}
 \end{equation*}
 from which \eqref{eq:directional_gradient_bound} follows with $C = \sqrt{2  d} $.
\end{proof*}

 First, the weak convergence of $\nabla_{\mathcal{D},\alpha} u$ will be studied, while successively a similar
 result will be obtained for $\nabla_{\mathcal{D},\alpha, \parallel} u$.
\begin{lemma}[Weak convergence of $\nabla_{\mathcal,\alpha} u$]
 Let $\Omega$ be a bounded open connected polyhedral subset of $\mathbb{R}^d$, $d
 \in \mathbb{N}^\star$. Let $\mathcal{D}$ be an admissible finite volume discretization in sense of Definition
 \eqref{def:polyhedral_mesh} and let $0 < \theta \leq \tilde\theta_\mathcal{D}$. Assume that there exists $u
 \in H_\mathcal{D}$ and a function $\underline{u} \in H_0^1(\Omega)$ such that $u$ tends to $\underline{u}$ in
 $L^2(\Omega)$ as $h_\mathcal{D} \rightarrow 0$, while $\|u\|_\mathcal{D}$ remains bounded. Then
 $\nabla_{\mathcal{D}} u$ weakly converges to $\alpha \nabla \underline{u}$ in $L^2(\Omega)^d$ as
 $h_\mathcal{D} \rightarrow 0$.
 \label{lem:gauss_gradient_weak_convergence}
\end{lemma}

\begin{proof*}
 Let $\varphi \in C_c^\infty(\Omega)$. Assume that $h_{\mathcal{D}}$ is small enough that, for all
 $K \in \mathcal{M}$ and $\bm{x} \in K$, if $\varphi(\bm{x}) \neq 0$ then $\mathcal{F}_{K,ext} = \varnothing$.
 Consider the term $T_1^\mathcal{D}$ defined as
 \begin{equation*}
  \begin{aligned}
   T_1^\mathcal{D} &= \int_\Omega P_\mathcal{D} \varphi(\bm{x}) \nabla_{\mathcal{D}} u(\bm{x}) \,\de\bm{x}
   = \sum_{K \in\ \mathcal{M}} |K| \,\varphi(\bm{x}_K) \, \nabla_K u \\
   &= \sum_{K|L \in \mathcal{F}_{int}} \left( \varphi(\bm{x}_K) \,d_{K,\sigma} + \varphi(\bm{x}_L) \,d_{L,\sigma} \right)
                                \alpha_{K|L} \, \tau_{K|L} \bm{n}_{K,\sigma} (u_L - u_K) .
  \end{aligned}
 \end{equation*}
 The first term between brackets can be rewritten as
 \begin{equation*}
  \begin{aligned}
   &\varphi(\bm{x}_K) \,d_{K,\sigma} + \varphi(\bm{x}_L) \,d_{L,\sigma} \\
   &\quad =   \vphantom{\frac{\bm{x}_K + \bm{x}_L}{2}}
        \varphi(\bm{x}_K) (\bm{x}_\sigma - \bm{x}_K) \cdot \bm{n}_{K,\sigma}    
       + \varphi(\bm{x}_L) (\bm{x}_L - \bm{x}_\sigma) \cdot \bm{n}_{K,\sigma} \\
   &\quad = \left( \frac{\varphi(\bm{x}_K) + \varphi(\bm{x}_L)}{2} (\bm{x}_L - \bm{x}_K) \right. \\
    &\hspace{3.0em}+\left. ( \varphi(\bm{x}_K) - \varphi(\bm{x}_L) ) \left( \bm{x}_\sigma - \frac{\bm{x}_K + \bm{x}_L}{2} \right)
      \right) \cdot \bm{n}_{K,\sigma}.
  \end{aligned}
 \end{equation*}
 The term $T_1^\mathcal{D}$ can be decomposed into a sum of two terms $T_1^\mathcal{D} = T_2^\mathcal{D} +
 T_3^\mathcal{D}$, where
 \begin{equation*}
  \begin{aligned}
   T_2^\mathcal{D} &= \sum_{K|L \in \mathcal{F}_{int}} \tau_{K|L}\, (u_L - u_K) \,\bm{n}_{K,\sigma}\, d_{K,L}\, 
                                              \frac{\varphi(\bm{x}_K) + \varphi(\bm{x}_L)}{2} , \\
   T_3^\mathcal{D} &= \sum_{K|L \in \mathcal{F}_{int}} \tau_{K|L}\, (u_L - u_K) \,\bm{n}_{K,\sigma}\,
                 \left( \bm{x}_\sigma - \frac{\bm{x}_K + \bm{x}_L}{2}\right) \cdot \bm{n}_{K,\sigma}.
  \end{aligned}
 \end{equation*}
 Starting with the analysis of term $T_3^\mathcal{D}$, by Cauchy-Schwartz and then triangle inequalities one
 gets
 \begin{equation*}
  \begin{aligned}
   &\left( T_3^\mathcal{D} \right)^2 \leq
      \sum_{K|L \in \mathcal{F}_{int}} \tau_{K|L} (u_L - u_K)^2\\
     &\quad \times \sum_{K|L \in \mathcal{F}_{int}} \tau_{K|L} \, (\varphi(\bm{x}_K) - \varphi(\bm{x}_L))^2
    \times \left| \left( \bm{x}_\sigma - \frac{\bm{x}_K + \bm{x}_L}{2} \right) \cdot \bm{n}_{K,\sigma} \right|^2
  \end{aligned}
 \end{equation*}
 in which, due to triangle inequality
 \begin{equation*}
  \left| \left( \bm{x}_\sigma - \frac{\bm{x}_K + \bm{x}_L}{2} \right) \cdot \bm{n}_{K,\sigma} \right|
  \leq \frac{1}{2} |\bm{x}_\sigma - \bm{x}_K| + \frac{1}{2} |\bm{x}_\sigma - \bm{x}_L|
  \leq h_{\mathcal{D}},
 \end{equation*}
 while due to mesh regularity
 \begin{multline*}
  \left| \bm{n}_{K,\sigma} - \frac{(\bm{x}_L - \bm{x}_K)}{d_{K,L}} \right|\\
  = \left| \bm{n}_{K,\sigma} - \frac{\bm{i}_{K,L}}{\bm{n}_{K,\sigma} \cdot \bm{i}_{K,L}} \right|
  \leq 1 + \left|\frac{\bm{i}_{K,L}}{\bm{n}_{K,\sigma} \cdot \bm{i}_{K,L}} \right|
  \leq 1 + \frac{1}{\theta},
 \end{multline*}
 from which it follows that
 \begin{equation*}
  \left( T_3^\mathcal{D} \right)^2 \leq
  C \, h_\mathcal{D} \,|\Omega|\, \|u\|_{\mathcal{D}}^2
 \end{equation*}
 with $C$ only depending on $d$, $\Omega$ and $\varphi$. Thus one concludes that
 $\lim_{h_\mathcal{D} \rightarrow 0} T_3^\mathcal{D} = 0$.
 Successively, compare $T_2^\mathcal{D}$ with the term
 \begin{equation*}
  \begin{aligned}
   T_4^\mathcal{D} &= -\int_\Omega u(\bm{x}) \nabla \varphi (\bm{x}) \,\de\bm{x} \\
   &= \sum_{K|L \in \mathcal{F}_{int}} (u_L - u_K) \int_{K|L} \varphi(\bm{x})
                                                               \, \bm{n}_{K,\sigma} \,\de\gamma(\bm{x}),
  \end{aligned}
 \end{equation*}
 which is such that
 \begin{equation*}
  \lim_{h_\mathcal{D} \rightarrow 0} T_4^\mathcal{D}
  = - \int_\Omega \underline{u}(\bm{x}) \nabla\varphi(\bm{x}) \,\de\bm{x}
  =  \int_\Omega \varphi(\bm{x}) \nabla \underline{u}(\bm{x}) \,\de\bm{x}.
 \end{equation*}
 Due to the fact that midpoint face interpolation is first order accurate
 \begin{equation*}
   \left| \frac{1}{|\sigma|} \int_{K|L} \varphi(\bm{x}) \,\de\gamma(\bm{x})
          - \frac{\varphi(\bm{x}_K) + \varphi(\bm{x}_L)}{2} \right|
   \leq h_{\mathcal{D}} \| \nabla \varphi \|_{L^\infty(\Omega)},
 \end{equation*}
 one has that
 \begin{equation*}
  \begin{aligned}
   &\left( T_4^\mathcal{D} - T_2^\mathcal{D} \right)^2\\
   &\hspace{1.0em}\leq \sum_{K|L \in \mathcal{F}_{int}} (\,|\sigma|\, \bm{n}_{K,\sigma})^2 (u_L - u_K)^2\\
   &\hspace{3.0em} \times \sum_{K|L \in \mathcal{F}_{int}} \left| \frac{1}{|\sigma|}
        \int_{K|L} \varphi(\bm{x}) \,\de\gamma(\bm{x}) - \frac{\varphi(\bm{x}_K) + \varphi(\bm{x}_L)}{2} \right|^2\\
   &\hspace{1.0em}\leq \sum_{K|L \in \mathcal{F}_{int}} |\sigma|^2 (u_L - u_K)^2
        \sum_{K|L \in \mathcal{F}_{int}}  h_{\mathcal{D}}^2 \| \nabla \varphi \|^2_{L^\infty(\Omega)},
  \end{aligned}
 \end{equation*}
 from which it follows that $\lim_{h_\mathcal{D} \rightarrow 0} \left( T_4^\mathcal{D} - T_2^\mathcal{D}
 \right)^2 = 0$.  Thus, $T_2^\mathcal{D}$ converges to $T_4^\mathcal{D}$ and, due to density of
 $C_c^\infty(\Omega)$ in $L^2(\Omega)$, $\nabla_{\mathcal{D}} u$ weakly converges to $\nabla
 \underline{u}$ as $h_{\mathcal{D}} \rightarrow 0$.
 Since
 \begin{equation*}
  \lim_{h_\mathcal{D} \rightarrow 0} T_4^\mathcal{D}
  = \int_\Omega \underline{u}(\bm{x}) \nabla\varphi(\bm{x}) \,\de\bm{x}
  = -\int_\Omega \varphi(\bm{x}) \nabla \underline{u}(\bm{x}) \,\de\bm{x},
 \end{equation*}
 by density of $C_c^\infty(\Omega)$ in $L^2(\Omega)$, one obtains the weak convergence of
 $\nabla_{\mathcal{D}} u (\bm{x})$ to $\nabla \underline{u} (\bm{x})$ as $h_{\mathcal{D}} \rightarrow 0$.
\end{proof*}

\begin{lemma}[Consistency of $\nabla_{\mathcal{D}}$]
 Let $\Omega$ be a bounded open connected polyhedral subset of $\mathbb{R}^d$, $d \in \mathbb{N}^\star$. Let
 $\mathcal{D}$ be an admissible finite volume discretization in sense of
 Definition \eqref{def:polyhedral_mesh} and let $0 < \theta \leq \tilde\theta_\mathcal{D}$. Let
 $\underline{u} \in C^2(\overline{\Omega})$ be such that $\overline{u} = 0$ on $\dd \Omega$. Then there exists
 $C$, depending only on $\Omega$, $\theta$, $\overline{u}$ and $\alpha$ such that
 \begin{equation}
   \| \nabla_{\mathcal{D}} P_{\mathcal{D}} \overline{u} - \nabla \overline{u} \|_{L^2(\Omega)^d}
   \leq C_3 h_{\mathcal{D}}
  \label{eq:gauss_gradient_consistency}
 \end{equation}
 \label{lem:gauss_gradient_consistency}
\end{lemma}

\begin{proof*}
 From Eq.\eqref{eq:discrete_gradient_gauss} for any $K \in \mathcal{M}$ one has
 \begin{equation*}
  \begin{aligned}
    |K| (\nabla_{\mathcal{D}} P_\mathcal{D} u)_K
   = &\sum_{L \in \mathcal{N}_K} \tau_{K|L}\, d_{K,\sigma}\, \bm{n}_{K,\sigma}\,
           (\underline{u}(\bm{x}_L)-\underline{u}(\bm{x}_K)) \\
     &\hspace{2.0cm} -\sum_{K \in \mathcal{M}} \sum_{\sigma \in \mathcal{F}_{K,ext}} \tau_{K,\sigma}\, d_{K,\sigma}\,
           \bm{n}_{K,\sigma}\, \underline{u}(\bm{x}_K)
  \end{aligned}
 \end{equation*}
 Let $(\nabla \underline{u})_K$ be the mean value of $\nabla \underline{u}$ over $K$
 \begin{equation*}
   (\nabla \underline{u})_K = \frac{1}{|K|} \int_K \nabla \underline{u}(\bm{x}) \,\de\bm{x}.
 \end{equation*}
 Due to the regularity of $\underline{u}$ and the homogeneous Dirichlet boundary conditions, the flux
 consistency error estimates include a constant $C$, only depending on $L^\infty$ norm of second derivatives
 of $\underline{u}$ (and of $\alpha$), such that for all $\sigma = K|L \in \mathcal{F}_{int}$, one has
 \begin{equation*}
   |e_\sigma| \leq C h_{\mathcal{D}} \quad
    \text{with} \quad e_\sigma = (\nabla \underline{u})_K \cdot {n}_{K,\sigma}
                              - \frac{\underline{u}(\bm{x}_L) - \underline{u}(\bm{x}_K)}{d_{K,L}}
 \end{equation*}
 while for all $\sigma \in \mathcal{F}_{ext}$ one has
 \begin{equation*}
   |e_\sigma| \leq C h_{\mathcal{D}} \quad
    \text{with} \quad e_\sigma = (\nabla \underline{u})_K \cdot {n}_{K,\sigma}
                              - \frac{- \underline{u}(\bm{x}_K)}{d_{K,\sigma}}
 \end{equation*}
 These flux consistency errors allow to recast $(\nabla_{\mathcal{D}} P_\mathcal{D} u)_K$ as
 \begin{equation*}
   |K| (\nabla_{\mathcal{D}} P_\mathcal{D} u)_K
    = \sum_{\sigma \in \mathcal{F}_{K,ext}} |\sigma| d_{K,\sigma} \bm{n}_{K,\sigma}
                                    (\nabla \underline{u})_K \cdot \bm{n}_{K,\sigma} + R_K 
 \end{equation*}
 where the consistency residual term is defined as
 \begin{equation*}
   R_K = - \sum_{\sigma \in \mathcal{F}_{K}} |\sigma| d_{K,\sigma} \bm{n}_{K,\sigma} e_\sigma.
 \end{equation*}
 From the geometrical identity valid for any vector $\bm{v} \in \mathbb{R}^d$ and for all $K \in \mathcal{M}$
 \begin{equation}
   \frac{1}{|K|} \sum_{\sigma \in \mathcal{F}_K} |\sigma| (\bm{x}_\sigma - \bm{x}_0) \bm{n}_{K,\sigma} \cdot \bm{v}
   = \bm{v} ,
 \end{equation}
 which is a direct consequence of the fact that each cell is a closed volume, it follows that
 \begin{equation*}
  |K| (\nabla_{\mathcal{D}} P_\mathcal{D} u)_K = |K| (\nabla \underline{u})_K + R_K
 \end{equation*}
 Due to flux consistency error estimates, it also follows that
 \begin{equation}
   |R_K| \leq C\, h_{\mathcal{D}} \sum_{\sigma \in \mathcal{F}_K} |\sigma|\, d_{K,\sigma}
          = C\, h_{\mathcal{D}}\, d\, |K|.
 \end{equation}
 As a consequence, one obtains that
 \begin{equation}
  \begin{aligned}
   &\sum_{K \in \mathcal{M}} |K| \left| (\nabla_{\mathcal{D}} P_{\mathcal{D}} u)_K - (\nabla \underline{u})_K \right|^2\\
   &\hspace{2.0em}\leq \sum_{K \in \mathcal{M}} |K|\, C^2 h_{\mathcal{D}}^2 d^2
    = |\Omega| C^2 h_{\mathcal{D}}^2 d^2.
  \end{aligned}
  \label{eq:density_weighted_gradient_consistency_a}
 \end{equation}

 Due to the regularity $\underline{u} \in C^2(\Omega)$, there exists another $C$, only dependent on $L^\infty$
 norm of the second derivatives of $\underline{u}$, such that
 \begin{equation}
  \sum_{K \in \mathcal{M}} \int_K \left| \nabla \underline{u} - (\nabla \underline{u})_K \right|^2
   \leq C h_\mathcal{D}^2.
  \label{eq:density_weighted_gradient_consistency_b}
 \end{equation}
 From Eqs.\eqref{eq:density_weighted_gradient_consistency_a}
 and \eqref{eq:density_weighted_gradient_consistency_b}, one gets the existence of $C_c$, only dependent on
 $\Omega$ and $\underline{u}$, such that \eqref{eq:gauss_gradient_consistency} holds.
\end{proof*}

\bibliographystyle{plain}
\bibliography{bib_fvm}

\end{document}